\documentclass{article}
\usepackage[utf8]{inputenc}
\usepackage{comment, amsmath, amssymb, amsthm, authblk, graphicx}
\usepackage[mathlines]{lineno}

\newtheorem{lem}{Lemma}[section]
\newtheorem{cor}{Corollary}[section]
\newtheorem{pro}{Proposition}[section]
\newtheorem{teo}{Theorem}[section]

\title{Expansive factors for geodesic flows of compact manifolds without conjugate points and with visibility universal covering}
\author[1]{Edhin F. Mamani\thanks{emamani@ufmg.br}}
\author[2]{Rafael Ruggiero\thanks{rorr@mat.puc-rio.br}}
\affil[1]{Instituto de Ciências Exatas ICEx, Universidade Federal de Minas Gerais, Av. Antônio Carlos 6627, Belo Horizonte 31270-901, Brazil.}
\affil[2]{Departamento de Matemática, Pontifícia Universidade Católica do Rio de Janeiro, Rua Marquês de São Vicente 225, Rio de Janeiro 22451-900, Brazil.}
\date{}
 
\begin{document}

\maketitle

\begin{abstract}
Let $(M,g)$ be a compact manifold without conjugate points and with visibility universal covering. We show that its geodesic flow has a time-preserving expansive factor which is topologically mixing and has a local product structure. As an application, assuming further the so-called entropy-gap we prove the uniqueness of the measure of maximal entropy for the geodesic flow. For the other results we restrict our setting assuming furthermore the continuity of Green bundles and the existence of a hyperbolic closed geodesic. In this new context, we deduce that Green bundles are uniquely integrable and are tangent to the smooth leaves of the horospherical foliations. Moreover, we prove the above expansive factor acts on a compact topological manifold and its geodesic flow has a unique measure of maximal entropy which has full support.
\end{abstract}

\section{Introduction}
The theory of geodesic flows of compact manifolds without conjugate points in higher dimension is a challenging topic in non-uniform dynamical systems. Diverse hypothesis are required to get important dynamical and ergodic properties of the geodesic flow. Some of these hypotheses has to do with the global geometry of the manifold, i.e., the geometry of the universal covering of the manifold. Its relevance was clear since Morse's work \cite{morse24} about the global geometry of geodesics of compact higher genus surfaces without conjugate points. In the 1970's, Eberlein and O'Neill \cite{eberoneil73} came out with the notion of visibility manifolds which allows the existence of regions of positive curvature, being thus a natural generalization of manifolds of negative curvature. Moreover, the universal covering of a compact higher genus surface without conjugate points is an important example of visibility manifold. So, it is natural to ask whether we can extend properties of the geodesic flow from the surface case to the higher dimensional visibility case. 

The article is divided in two parts in the context of compact manifolds without conjugate points and with visibility universal covering: firstly we deal with the uniqueness of the measure of maximal entropy for geodesic flows. Secondly, assuming further the continuity of Green bundles and the existence of a hyperbolic closed geodesic, we study the uniqueness, smoothness and tangency of the horospherical foliations. 

The uniqueness of the measure of maximal entropy for geodesic flows of compact manifolds of negative curvature was showed by Bowen in 1972 \cite{bowen73maxi}. This result was extended to compact rank-1 manifolds of non-positive curvature by Knieper in 1998 \cite{knie98}. In 2018, this property was proved for compact higher genus surfaces without focal points by Gelfert and Ruggiero using an expansive factor of the geodesic flow \cite{gelf19}. In 2020, Climenhaga, Knieper and War showed the same conclusion for certain family of compact manifolds without conjugate points using Climenhaga-Thompson criterion \cite{clim16}. In this work, we follow Gelfert-Ruggiero's method since it is a shorter path to show the uniqueness feature. Our first contribution is an extension of results founded originally for compact surfaces without focal points.
\begin{teo}\label{i1}
Let $M$ be a compact $C^{\infty}$ $n$-dimensional manifold without conjugate points and with visibility universal covering. Then, the geodesic flow is time-preserving semi-conjugate to a continuous expansive flow $\psi_t$ acting on a compact metric space of topological dimension at least $n-1$. Moreover, $\psi_t$ is topologically mixing and has a local product structure.
\end{teo}
In order to show the uniqueness property in our context, we assume the so-called entropy-gap. This assumptions came from Climenhaga-Thompson's work \cite{clim16} and it was also used in Climenhaga-Knieper-War's Theorem. Our second contribution follows.
\begin{teo}\label{i2}
Let $M$ be a compact $C^{\infty}$ $n$-manifold without conjugate points which has a visibility universal covering and $\phi_t$ be its geodesic flow. If $\sup \{ h_{\mu}(\phi_1): \mu \text{ is supported on } T_1M\setminus \mathcal{R}_0 \}<h(\phi_1)$ then $\phi_t$ has a unique measure of maximal entropy.
\end{teo}
We highlight that our setting is more general than Climenhaga-Knieper-War context because we eliminate and generalize some hypothesis of their theorem. For more details see Section \ref{ent}.

By Green's work we know that Green bundles always exist for compact manifolds without conjugate points \cite{green58}. These bundles are measurable, Lagrangian and invariant by the geodesic flow. In the case of Anosov geodesic flows, Green bundles are the dynamical invariant bundles of Anosov dynamics. The continuity of Green bundles holds in several settings: manifolds of non-positive curvature, manifolds without focal points and manifolds of bounded asymptote.

In the hyperbolic case, i.e., when the manifold has negative curvature, the horospherical foliations are the only continuous invariant foliations with smooth leaves tangent to Green bundles which are just the invariant bundles given by Anosov's theory \cite{anos67}. This property was extended to compact surfaces of non-positive curvature by Eberlein \cite{heint77}. However, when some regions of positive curvature are allowed, the smoothness of horospherical leaves is a longstanding open problem. Moreover, Knieper \cite{knip86} observed that when Green bundles are continuous, they are integrable, but neither necessarily tangent to horospherical leaves, nor uniquely integrable. In the case of compact higher genus surfaces without conjugate points, it is not known neither if Green bundles are tangent to the horospherical foliations, nor if horospherical leaves are smooth. 

On the other hand, Barbosa and Ruggiero \cite{barbo07} showed that horospherical foliations of compact higher genus surfaces without conjugate points are the only continuous foliations of the unit tangent bundle invariant by the geodesic flow. As far as we know, this is the state of the art of the problem in the theory of manifolds without conjugate points admitting regions of positive curvature. Our following statement gives a partial answer to these questions in our higher dimensional setting. In particular, we extend Barbosa-Ruggiero's Theorem \cite{barbo07}, Rosas-Ruggiero's Theorem \cite{rosas03} and some results in Gelfert-Ruggiero's work \cite{gelf20} to higher dimensions.

\begin{teo}\label{mago}
Let $M$ be a compact $C^{\infty}$ $n$-manifold without conjugate points, with visibility universal covering and with continuous Green bundles. Suppose the geodesic flow has a periodic hyperbolic point $\theta\in T_1M$ then
\begin{enumerate}
\item The set of points with non-zero Lyapunov exponents in directions transverse to the flow agrees almost everywhere with an open dense set, with respect to Liouville measure.
\item Hyperbolic periodic points are dense on $T_1M$.
\item Green bundles are uniquely integrable and tangent to the smooth horospherical foliations $\mathcal{F}^s$ and $\mathcal{F}^u$.
\end{enumerate}
\end{teo}
Clearly the setting of Theorem \ref{mago} is a particular case of Theorems \ref{i1} and \ref{i2}. So, these theorems provide the same conclusions in the second part of the article. We can say a bit more about the properties of the objects involved in the above theorems.
\begin{teo}\label{magi}
Under the hypothesis of Theorem \ref{mago} we have
\begin{enumerate}
\item The geodesic flow $\phi_t$ is time-preserving semi-conjugate to a continuous expansive flow acting on a compact topological $(2n-1)$-manifold $X$ which is topologically mixing and has a local product structure.
\item $\phi_t$ has a unique measure of maximal entropy which has full support.
\end{enumerate}
\end{teo}
The organization of the paper is as follows. Section 2 states notations and recalls some preliminary concepts and results. In Section 3, we recall the construction of the factor flow and prove some basic properties. In Section 4, we construct the basis of neighborhoods in our higher dimensional context. Section 5 sketches the proof of the dynamical and ergodic properties of the factor flow. In Section 6, we show the uniqueness of the measure of maximal entropy. Section 7 begins the second part of the article which deals with the manifolds considered in Theorem \ref{mago} and proves items (1) and (2) of theorem \ref{mago}. Section 8 addresses the regularity and tangency problem of the horospherical foliations and shows item (3) of Theorem \ref{mago}. In section 9, we study the regularity of the quotient space and prove item (1) of Theorem \ref{magi}. Section 10 shows that manifolds in Theorem \ref{magi} have geodesic flows with unique measures of maximal entropy. 
    
\section{Preliminaries}
\subsection{Compact manifolds without conjugate points}\label{m}
We introduce the general background and main notations we will use throughout the paper. Let $(M,g)$ be a $C^{\infty}$ compact connected Riemannian manifold, $TM$ be its tangent bundle and $T_1M$ be its unit tangent bundle. Consider the universal covering $\tilde{M}$ of $M$, the covering map $\pi:\tilde{M}\to M$ and the natural projection $d\pi:T\tilde{M}\to TM$. The universal covering $(\tilde{M},\tilde{g})$ is a complete Riemannian manifold with the pullback metric $\tilde{g}$. A manifold $M$ has no conjugate points if the exponential map $\exp_p$ is non-singular at every $p\in M$. In particular, $\exp_p$ is a covering map for every $p\in M$ (p. 151 of \cite{doca92}).

We denote by $\nabla$ the Levi-Civita connection of $(M,g)$. A geodesic is a smooth curve $\gamma\subset M$ with $\nabla_{\dot{\gamma}}\dot{\gamma}=0$. For every $\theta=(p,v)\in TM$, $\gamma_{\theta}$ is the unique geodesic satisfying $\gamma_{\theta}(0)=p$ and $\dot{\gamma}_{\theta}(0)=v$. The geodesic flow $\phi_t$ is defined by 
$$\phi: \mathbb{R}\times TM\to TM, \qquad (t,\theta)\mapsto \phi_t(\theta)=\dot{\gamma}_{\theta}(t).$$
If every geodesic is parametrized by arc-length, we can restrict the geodesic flow to $T_1M$. 

We now define a Riemannian metric on the tangent bundle $TM$ (Section 1.3 of \cite{pater97}). Denote by $P:TM\to M$ and $\tilde{P}:T\tilde{M}\to \tilde{M}$ the corresponding canonical projections. For every $\theta=(p,v)\in TM$, the Levi-Civita connection induces the so-called connection map $C_{\theta}:T_{\theta}TM\to T_pM$. These linear maps provide the linear isomorphism $T_{\theta}TM\to T_pM\oplus T_pM$ with $\xi\mapsto (d_{\theta}P(\xi),C_{\theta}(\xi))$. We define the horizontal subspace by $\mathcal{H}(\theta)=\ker(C_{\theta})$ and the vertical subspace by $\mathcal{V}(\theta)=\ker(d_{\theta}P)$. These subspaces decompose the tangent space by $T_{\theta}TM=\mathcal{H}(\theta)\oplus \mathcal{V}(\theta)$. For every $\xi,\eta\in T_{\theta}TM$, the Sasaki metric is defined by
\begin{equation}\label{sasa}
    \langle \xi,\eta \rangle_s = \langle d_{\theta}P(\xi), d_{\theta}P(\eta) \rangle_p + \langle C_{\theta}(\xi), C_{\theta}(\eta) \rangle_p.
\end{equation}
This metric induces a Riemannian distance $d_s$ usually called Sasaki distance.

For every $\theta\in T_1M$, denote by $G(\theta)\subset T_{\theta}T_1M$ the subspace tangent to the geodesic flow at $\theta$. Let $N(\theta)\subset T_{\theta}T_1M$ be the subspace orthogonal to $G(\theta)$ with respect to the Sasaki metric. For every $\theta\in T_1M$, $H(\theta)=\mathcal{H}(\theta)\cap N(\theta)$ and $V(\theta)=\mathcal{V}(\theta)\cap N(\theta)$.
From the above decomposition we have
\[  T_{\theta}T_1M=H(\theta)\oplus V(\theta)\oplus G(\theta) \quad \text{ and }\quad N(\theta)=H(\theta)\oplus V(\theta). \]
So, every $\xi\in N(\theta)$ has decomposition $\xi=(\xi_h,\xi_v)\in H(\theta)\oplus V(\theta)$. We call $\xi_h$ and $\xi_v$ the horizontal and vertical components of $\xi$ respectively.

\subsection{Horospheres}\label{h}
We assume that $(M,g)$ is a $C^{\infty}$ compact n-manifold without conjugate points. Let us introduce important asymptotic objects in the universal covering (\cite{esch77} and part II of \cite{pesin77}). Let $\theta\in T_1\tilde{M}$ and $\gamma_{\theta}$ be the geodesic induced by $\theta$. Define the forward Busemann function by
$$b_{\theta}: \tilde{M}\to \mathbb{R}, \qquad p\mapsto b_{\theta}(p)=\lim_{t\to \infty}d(p,\gamma_{\theta}(t))-t.$$
From now on, for every $\theta=(p,v)\in T_1\Tilde{M}$ we denote $-\theta:=(p,-v)\in T_1\Tilde{M}$. The stable and unstable horosphere of $\theta$ are defined by
\[  H^+(\theta)=b_{\theta}^{-1}(0)\subset \tilde{M} \quad \text{ and }\quad  H^-(\theta)=b_{-\theta}^{-1}(0)\subset \tilde{M}. \]
We lift these horospheres to $T_1\tilde{M}$. Denote by $\nabla b_{\theta}$ the gradient vector field of $b_{\theta}$ and define the sets
\[  \mathcal{\tilde{F}}^s(\theta)=\{ (p,-\nabla_pb_{\theta}): p\in H^+(\theta) \}\quad \text{ and }\quad \mathcal{\tilde{F}}^u(\theta)= \{ (p,\nabla_pb_{-\theta}): p\in H^-(\theta) \}. \]
These sets project onto the horospheres by the canonical projection $\tilde{P}$. For every $\theta\in T_1\tilde{M}$, define the stable and unstable families of sets of $\theta$ by
$$\mathcal{\tilde{F}}^s=(\mathcal{\tilde{F}}^s(\theta) )_{\theta\in T_1\tilde{M}} \quad \text{ and }\quad \mathcal{\tilde{F}}^u=( \mathcal{\tilde{F}}^u(\theta) )_{\theta\in T_1\tilde{M}}.$$
We also define the center stable and center unstable sets of $\theta$ by 
$$\mathcal{\tilde{F}}^{cs}(\theta)=\bigcup_{t\in \mathbb{R}} \mathcal{\tilde{F}}^s(\phi_t(\theta))  \quad \text{ and }\quad \mathcal{\tilde{F}}^{cu}(\theta)=\bigcup_{t\in \mathbb{R}} \mathcal{\tilde{F}}^u(\phi_t(\theta)).$$

The integral flow $\sigma_t^{\theta}:\tilde{M}\to \tilde{M}$ of $-\nabla b_{\theta}$ is called the Busemann flow and its integral curves are called Busemann asymptotes of $\gamma_{\theta}$. In particular, $\gamma_{\theta}$ is an integral curve of $\sigma_t^{\theta}$. In contrast, two curves $\gamma_1,\gamma_2\subset \tilde{M}$ are asymptotic if $d(\gamma_1(t),\gamma_2(t))\leq C$ for every $t\geq 0$ and some $C>0$. If in addition $\gamma_1(-t)$ and $\gamma_2(-t)$ are asymptotic then $\gamma_1$ and $\gamma_2$ are called bi-asymptotic. Let us state some basic properties of these objects.

\begin{pro}[\cite{esch77,pesin77}]\label{horo}
Let $M$ be a compact manifold without conjugate points. Then, for every $\theta\in T_1\tilde{M}$,
\begin{enumerate}
    \item Busemann functions $b_{\theta}$ are $C^{1,L}$ with $L$-Lipschitz unitary gradient for a uniform constant $L>0$ \cite{knip86}.
    \item Horospheres $H^+(\theta),H^-(\theta)\subset \tilde{M}$ and sets $\mathcal{\tilde{F}}^s(\theta),\mathcal{\tilde{F}}^u(\theta)\subset T_1\tilde{M}$ are Lipschitz continuous embedded $n-1$-submanifolds.
    \item Busemann asymptotes to $\gamma_{\theta}$ are always orthogonal to $H^+(\theta)$. 
    \item Horospheres are equidistant: for every $t,s\in\mathbf{R}$ and $q\in H^+(\gamma_{\theta}(t))$, 
    $$d(q,H^+(\gamma_{\theta}(s)))=|t-s|.$$
\end{enumerate}
\end{pro}

We can define the above objects in the case of $T_1M$. For every $\theta\in T_1M$, 
$$\mathcal{F}^*(\theta)=d\pi (\mathcal{\tilde{F}}^*(\tilde{\theta}))\subset T_1M \quad \text{ and }\quad \mathcal{F}^*=d\pi (\mathcal{\tilde{F}}^*), \quad *=s,u,cs,cu;$$
for any lift $\tilde{\theta}\in T_1\tilde{M}$ of $\theta$. 

\subsection{Visibility manifolds}\label{v}
In this subsection we introduce some dynamical and geometric properties of visibility manifolds. Let $M$ be a simply connected Riemannian manifold without conjugate points. For every $x,y,z\in M$, denote by $[x,y]$ the unique geodesic segment joining $x$ to $y$ and by $\sphericalangle_z(x,y)$ the angle at $z$ formed by $[z,x]$ and $[z,y]$. We say that $M$ is a visibility manifold if for every $z\in M$ and every $\epsilon>0$ there exists $R(\epsilon, z)>0$ such that 
$$\text{ if }x,y\in M \text{ with }d(z,[x,y])>R(\epsilon, z) \quad \text{ then }\quad \sphericalangle_z(x,y)<\epsilon.$$
If $R(\epsilon,z)$ does not depend on $z$ then $M$ is called a uniform visibility manifold. We note that every compact manifold without conjugate points has a uniform visibility universal cover \cite{eber72}. 

In 1973, Eberlein extended some transitivity properties of the geodesic flow to the setting of compact manifolds without conjugate points. A foliation is called minimal if its leaves are dense.

\begin{teo}[\cite{eber72,eber73neg2}]\label{v1}
Let $M$ be a compact manifold without conjugate points and with visibility universal covering $\Tilde{M}$. Then
\begin{enumerate}
    \item The families $\mathcal{F}^s$ and $\mathcal{F}^u$ are continuous minimal foliations of $T_1M$ invariant by the geodesic flow: for every $t\in \mathbb{R}$,
    $$\phi_t(\mathcal{F}^*(\theta))=\mathcal{F}^*(\phi_t(\theta)), \qquad *=s,u.$$
    \item The geodesic flow $\phi_t$ is topologically mixing.
    \item There exist $A,B>0$ such that for every $\theta\in T_1\tilde{M}$ and every $\eta\in \mathcal{\tilde{F}}^s(\theta)$, 
    $$d_s(\phi_t(\theta),\phi_t(\eta))\leq Ad_s(\theta,\eta)+B \quad \text{ for every } t\geq 0.$$
\end{enumerate}
\end{teo}
Item 1 clearly extends to the setting $T_1\tilde{M}$. The families $\mathcal{\tilde{F}}^s$ and $\mathcal{\tilde{F}}^u$ are called the stable and unstable horospherical foliations of $T_1\tilde{M}$. Thus, for every $\theta\in T_1\tilde{M}$, $\mathcal{\tilde{F}}^s(\theta)$ and $\mathcal{\tilde{F}}^u(\theta)$ are called the stable and unstable horospherical leaf of $\theta$. Item 3 says that horospherical foliations still have some kind of weak hyperbolicity. 

We now deal with intersections between horospherical leaves. For every $\theta\in T_1M$, denote by
$I(\theta)=H^+(\theta)\cap H^-(\theta)$ and by $\mathcal{I}(\theta)=\mathcal{F}^s(\theta)\cap \mathcal{F}^u(\theta)$. The canonical projection $\tilde{P}$ maps $\mathcal{I}(\theta)$ onto $I(\theta)$. We say that $\theta\in T_1M$ is expansive if $\mathcal{I}(\theta)=\{ \theta\}$ otherwise $\theta$ is non-expansive. The expansive set is defined by
$$\mathcal{R}_0=\{ \theta\in T_1M: \mathcal{I}(\theta)=\{ \theta\} \}.$$
The complement $T_1M\setminus\mathcal{R}_0$ is called the non-expansive set. These sets help to characterize the dynamical and ergodic behaviour of the geodesic flow. 

For intersections of stable and unstable horospherical leaves of different points, the following property holds for visibility manifolds.

\begin{teo}[\cite{eber72}]\label{conec}
If $M$ is a compact manifold without conjugate points and with visibility universal covering $\Tilde{M}$ then for every $\theta,\xi \in T_1\tilde{M}$ with $\theta\not\in \mathcal{\tilde{F}}^{cu}(\xi)$ there exists $\eta_1,\eta_2\in T_1\tilde{M}$ such that
$$\mathcal{\tilde{F}}^s(\theta)\cap \mathcal{\tilde{F}}^{cu}(\xi)=\mathcal{\tilde{I}}(\eta_1) \quad \text{ and }\quad  \mathcal{\tilde{F}}^s(\xi)\cap \mathcal{\tilde{F}}^{cu}(\theta)=\mathcal{\Tilde{I}}(\eta_2).$$
\end{teo}
This formula can be rephrased in terms of unstable horospherical leaves and central stable sets: there exist $t_1,t_2\in \mathbb{R}$ such that
$$\mathcal{\tilde{F}}^{cs}(\theta)\cap \mathcal{\tilde{F}}^u(\xi)=\mathcal{\Tilde{I}}(\Tilde{\phi}_{t_1}(\eta_1)) \quad \text{ and }\quad \mathcal{\tilde{F}}^{cs}(\xi)\cap\mathcal{\tilde{F}}^u(\theta)=\mathcal{\Tilde{I}}(\Tilde{\phi}_{t_2}(\eta_2)).$$
The above intersections are called the heteroclinic connections of the geodesic flow. 

\subsection{Jacobi fields and Green bundles}

Let $\theta\in T_1M$ and $\gamma_{\theta}$ be the geodesic induced by $\theta$. A vector field $J$ along $\gamma_{\theta}$ is called a Jacobi field if it satisfies the Jacobi equation
$$J''(t)+R(\dot{\gamma}_{\theta}(t),J(t))\dot{\gamma}_{\theta}(t)=0,$$
where $R$ is the curvature tensor induced by the metric $g$. 

From this equation, we see that a Jacobi field is completely determined by its initial conditions. Moreover, a direct calculation shows that $J(t)$ is orthogonal to $\dot{\gamma}_{\theta}(t)$ for all $t\in \mathbb{R}$ if and only if $J(0)$ and $J'(0)$ both are orthogonal to $\dot{\gamma}_{\theta}(0)=\theta$. In this case, $J$ is called an orthogonal Jacobi field. We denote by $\mathcal{J}_{\theta}$ the $(2n-2)$-dimensional vector subspace of orthogonal Jacobi fields on $\gamma_{\theta}$.

From subsection \ref{m}, the orthogonal decomposition allows to write any vector $\xi\in T_{\theta}T_1M$ as $\xi=(\xi_h,\xi_v)$. Let $J_{\xi}$ be the Jacobi field on $\gamma_{\theta}$ determined by the initial conditions $J_{\xi}(0)=\xi_h$ and $J'_{\xi}(0)=\xi_v$. By above, we see that $J_{\xi}\in \mathcal{J}_{\theta}$ if and only if $\xi\in N(\theta)$. This relation can be extended to an isomorphism.
\begin{pro}\label{j1}
Let $M$ be a compact manifold without conjugate points. Then, for every $\theta=(p,v)\in T_1M$ we have
\begin{enumerate}
    \item The map $\xi\mapsto J_{\xi}$ is an isomorphism between $N(\theta)$ and $\mathcal{J}_{\theta}$ \cite{eber73.1}.
    \item For every $\xi \in N(\theta)$ and every $t\in \mathbb{R}$ \cite{pater12},
    $$d_{\theta}\phi_t(\xi)=(J_{\xi}(t),J'_{\xi}(t)) \quad \text{ hence }\quad \|d_{\theta}\phi_t(\xi)\|^2_v=\|J_{\xi}(t)\|^2_p+\|J'_{\xi}(t)\|^2_p.$$    
\end{enumerate}
\end{pro}
In 1958 Green \cite{green58} introduced a distinguished class of Jacobi fields that always exist for compact manifolds without conjugate points. Let $\theta=(p,v)\in T_1M$, $\gamma_{\theta}$ be the geodesic induced by $\theta$. For every $T\in \mathbb{R}$ and every unit $w\in v^{\perp}\subset T_pM$, consider the Jacobi field $J_T\in \mathcal{J}_{\theta}$ with boundary conditions $J_T(0)=w$ and $J_T(T)=0$. Green showed that limit Jacobi fields always exist when $T\to \pm\infty$. We call to
$$J_w^s=\lim_{T\to \infty}J_T \quad \text{ and } \quad J_w^u=\lim_{T\to -\infty}J_T$$
the stable and unstable Jacobi fields on $\gamma_{\theta}$ with initial condition $J_w^s(0)=J_w^u(0)=w$. These fields never vanish and form two vector subspaces of $\mathcal{J}_{\theta}$. Proposition \ref{j1} allows to lift the vector subspaces of stable and unstable Jacobi fields on $\gamma_{\theta}$ to vector subspaces $G^s(\theta)$ and $G^u(\theta)$ of $T_{\theta}T_1M$. The collections $G^s=\bigcup_{\theta\in T_1M}G^s(\theta)$ and $G^u=\bigcup_{\theta\in T_1M}G^u(\theta)$ are $(n-1)$-dimensional sub-bundles of the tangent bundle of $T_1M$ which are called the stable and unstable Green bundles respectively. Green bundles are measurable and invariant by the derivative of the geodesic flow $d\phi_t$ \cite{eber73.1}. We say that Green bundles are continuous if $G^s(\theta)$ and $G^u(\theta)$ depend continuously on $\theta\in T_1M$. 

In connection with set $\mathcal{R}_0$ from subsection \ref{v}, we define the important set
$$\mathcal{R}_1=\{ \theta\in T_1M: G^s(\theta)\cap G^u(\theta)= \{0\}\},$$
i.e., the set where Green bundles are linearly independent. We have $\mathcal{R}_1\subset \mathcal{R}_0$ in particular contexts: non-positive curvature manifolds and manifolds without focal points. In the general case the inclusion remains as an important open problem.

For vectors belonging to the stable and unstable Green bundles we have the following improvement of Proposition \ref{j1}(2).
\begin{pro}[Proposition 2.11 of \cite{eber73.1}]
Let $M$ be a compact manifold without conjugate points. Then, there exists $K>0$ such that for every $\xi\in G^s(G^u)$ and for every $t\geq 0 (t\leq 0)$, $\|J_{\xi}(t)\|\leq \|d\phi_t(\xi)\|\leq K \|J_{\xi}(t)\|$.
\end{pro}

\subsection{Generalized flat strip theorem}

In this subsection we give more properties of the intersections $\mathcal{I}(\theta)$. We first give some definitions of global geometry properties. Let $M$ be a simply connected manifold without conjugate points. We say that geodesics rays diverge in $M$ if for every $p\in M$ and every $\epsilon,C>0$ there exists $T(p,\epsilon,C)$ such that for every two geodesics rays $\gamma,\beta \subset M$ with same starting point $p$ and $\sphericalangle_p(\dot{\gamma}(0),\dot{\beta}(0))\geq \epsilon$ we have $d(\gamma(t),\beta(t))\geq C$ for every $t\geq T(p,\epsilon,C)$. The geodesics diverge uniformly if $T(p,\epsilon,C)$ does not depend on $p$. We say that $M$ is \textbf{quasi-convex} if there exist constants $A, B>0$ such that for every two geodesic segments $\gamma : [t_1, t_2] \to \tilde{M},\beta : [s_1,s_2] \to \tilde{M}$ it holds that
$$ d_H(\gamma , \beta ) \leq A \sup \{ d(\gamma(t_1), \beta(s_1)), d(\gamma(t_2), \beta(s_2))\}+ B, $$
where $d_H$ is the Hausdorff distance. Quasi-convexity with uniform divergence of geodesic rays provide a general framework for establishing some asymptotic geodesic properties.

\begin{teo}[\cite{riff18}]\label{rif}
Let $M$ be a compact manifold without conjugate points and $\tilde{M}$ be its universal covering. Suppose that $\tilde{M}$ is quasi-convex and geodesic rays diverge. Then, for every $\theta,\eta\in T_1\tilde{M}$ with bi-asymptotic geodesics $\gamma_{\theta},\gamma_{\eta}$ there exists a connected set $\Sigma\subset I(\theta)$ containing the geodesic starting points such that any geodesic $\beta$ with $\beta(0)=q\in \Sigma$ and $\dot{\beta}(0)=-\nabla b_q^{\theta}$ is bi-asymptotic to both $\gamma_{\theta},\gamma_{\eta}$. In particular, the set 
$$S=\bigcup_{t\in\mathbb{R},q\in \Sigma}\gamma_{(q,-\nabla b_q^{\theta})}(t)$$
is homeomorphic to $\Sigma \times \mathbb{R}$.
\end{teo}

The set $S$ is the topological analogous of the flat strip in the flat strip Theorem \cite{esch77, pesin77}. Roughly, the theorem says that every two bi-asymptotic geodesics bound a connected set (the strip) composed of bi-asymptotic geodesics. We can rephrase the conclusion in terms of horospherical leaves: there exists a connected set $\tilde{\Sigma}\subset \mathcal{I}(\theta)$ such that any $\gamma_{\xi}$ with $\xi\in \tilde{\Sigma}$ is bi-asymptotic to $\gamma_{\theta}$. The following theorem gives the connection with our context.

\begin{teo}\label{unif}
Let $M$ be a compact manifold without conjugate points and with visibility universal covering $\Tilde{M}$. Then,
    \begin{enumerate}
        \item $\tilde{M}$ is quasi-convex \cite{eber72}.
        \item Geodesic rays diverge uniformly in $\tilde{M}$ \cite{rugg94,rugg07}.
        \item There exists a uniform constant $\Tilde{Q}>0$ such that for every two bi-asymptotic geodesics $\beta,\gamma\subset \tilde{M}$ it holds that $d_H(\beta, \gamma)\leq \Tilde{Q}$ \cite{eber72}.
    \end{enumerate}
\end{teo}
This theorem gives some information about the global geometry of the generalized strips.
\begin{cor}[\cite{riff18}]\label{morse}
Let $M$ be a compact manifold without conjugate points, with visibility universal covering $\Tilde{M}$ and $\Tilde{Q}>0$ from Theorem \ref{unif}. Then there exists an uniform constant $Q\leq \Tilde{Q}$ such that for every $\eta \in T_1\Tilde{M}$,
the sets $I(\theta),\mathcal{\Tilde{I}}(\eta)$ are compact connected sets with $Diam(I(\theta))\leq Q$ and $Diam(\mathcal{\Tilde{I}}(\theta))\leq \Tilde{Q}$. 
\end{cor}
The uniform constant $\Tilde{Q}$ is usually called Morse's constant and depends only on the manifold $M$. Theorem \ref{rif} also gives a criterion to determine whether a vector belongs to an intersection $\mathcal{\Tilde{I}}(\theta)$ in terms of bi-asymptoticity.
\begin{cor}[\cite{riff18}]\label{key}
Let $M$ be a compact manifold without conjugate points and with visibility universal covering $\Tilde{M}$. For every $\theta\in T_1\tilde{M}$, if $\eta=(q,w)\in \mathcal{\tilde{F}}^s(\theta)$ and $\gamma_{\eta}$ is bi-asymptotic to $\gamma_{\theta}$ then 
\[ \eta \in \mathcal{\Tilde{I}}(\theta)=\mathcal{\tilde{F}}^s(\theta)\cap\mathcal{\tilde{F}}^u(\theta) \quad \text{ and }\quad q\in I(\theta)=H^+(\theta)\cap H^-(\theta). \]
\end{cor}

\subsection{Continuous parametrizations of horospherical leaves by polar coordinates}\label{parame}

In this subsection we will give polar coordinates to pieces of horospheres. These pieces will be large enough to contain certain non-trivial intersections $\mathcal{I}(\theta)$. We first recall the approximation of horospheres by geodesic spheres.
\begin{teo}[\cite{rugg03}]\label{T}
Let $M$ be a simply connected manifold where geodesic rays diverge uniformly. Then, for every $R,\epsilon>0$ there exists $T(R,\epsilon)>0$ such that for every $\theta=(p,v)\in T_1\tilde{M}$
$$d_H(H^+(\theta)\cap \Bar{B}_R(p),S_T(\gamma_{\theta}(T))\cap \Bar{B}_R(p))\leq \epsilon,$$
where $S_T(\gamma_{\theta}(T))$ is the geodesic sphere of radius $T$ centered at $\gamma_{\theta}(T)$ and $\Bar{B}_R(p)$ is the closed ball of radius $R$ centered at $p$.
\end{teo}
Let $Q>0$ be the uniform constant given by Corollary \ref{morse}. Choosing $\epsilon_0>0$ small enough, Theorem \ref{T} gives a time $T(2Q,\epsilon_0)>0$ such that for every $\theta=(p,v)\in T_1\tilde{M}$ 
$$d_H(H^+(\theta)\cap \Bar{B}_{2Q}(p),S_T(\gamma_{\theta}(T))\cap \Bar{B}_{2Q}(p))\leq \epsilon_0.$$
In particular, by Corollary \ref{T} we see that for every $\theta\in T_1\Tilde{M}$, $I(\theta)\subset H^+(\theta)\cap \Bar{B}_{2Q}(p)$. 

\begin{figure}[ht]
    \centering
    \includegraphics[width=\textwidth]{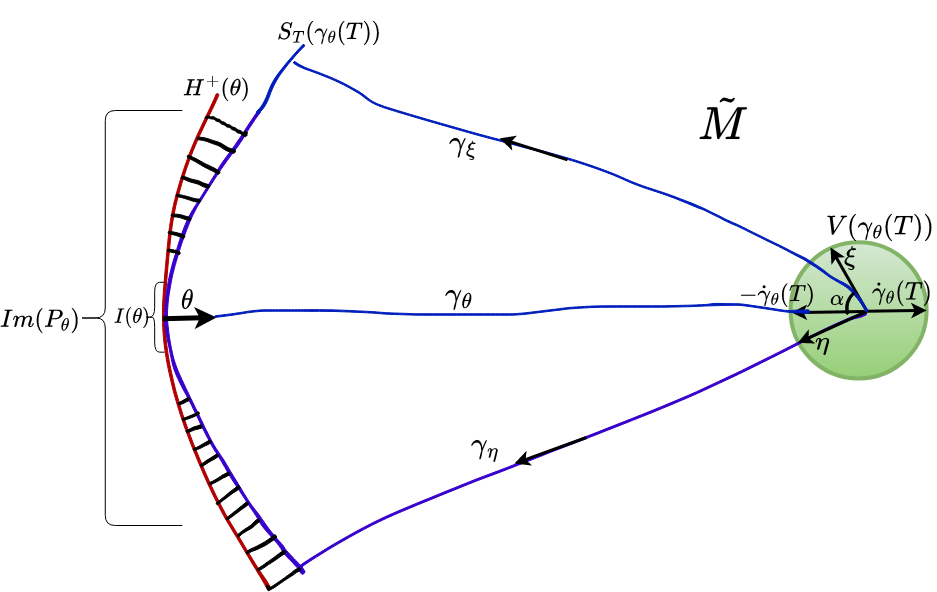}
    \caption{Construction of the parametrization of the stable horosphere of $\theta$.}
    \label{esfe}
\end{figure}

Using $T=T(2Q,\epsilon_0)$ we will define a projection map. For every $\theta\in T_1\Tilde{M}$, choose an orthonormal basis $\mathcal{B}=\{e_1,\ldots,e_n\}$ for $T_{\gamma_{\theta}(T)}\Tilde{M}$ such that $e_1=-\dot{\gamma}_{\theta}(T)$. Recall that for every $\theta=(p,v)\in T_1\Tilde{M}$, its vertical fiber is defined by $V(\theta)=\{ v\in T_p\Tilde{M}: |v|=1 \}$. It is well-known that $\mathcal{B}$ provides a parametrization of $V(\gamma_{\theta}(T))\setminus \{ \dot{\gamma}_{\theta}(T) \}$ by polar or exponential coordinates. This parametrization consists of $(n-1)$ oriented angles measured from $e_1$ to $e_i$ axis for $i=2,\dots, n$. We can see all the construction for the surface case in Figure \ref{esfe}. These angles also parametrize $H^+(\theta)\cap \Bar{B}_{2Q}(p)$ and $S_T(\gamma_{\theta}(T))\cap \Bar{B}_{2Q}(p)$ as we will see. Choosing some $a<\pi$, we define a projection map
$$P_{\theta}: [0,a]^{n-1}\longrightarrow H^+(\theta), \quad (\alpha_1,\ldots, \alpha_{n-1})\mapsto P_{\theta}(\alpha_1,\ldots, \alpha_{n-1})$$
as follows. Let $\eta\in V(\gamma_{\theta}(T))$ be the unique vector making angle $\alpha_i$ with $e_i$-axis for $i=1,\ldots, n-1$ and so $P_{\theta}(\alpha_1,\ldots, \alpha_{n-1})=\gamma_{\eta}\cap H^+(\theta)$. We choose $a>0$ so that $Im(P_{\theta})\subset \Bar{B}_{2Q}(p)$. Note that geodesics at time $T$ generated in this way form part of a geodesic sphere:
$$\{ \gamma_{\eta}(T): \eta \text{ induced by  angles } (\alpha_1,\ldots ,\alpha_{n-1}) \in [0,a]^{n-1}  \} \subset S_T(\gamma_{\theta}(T)).$$
Moreover, Theorem \ref{T} says that the geodesic sphere $S_T(\gamma_{\theta}(T))\cap \Bar{B}_{2Q}(p)$ is $\epsilon$-close to $H^+(\theta)\cap \Bar{B}_{2Q}(p)$ in the Hausdorff distance. This geodesic sphere is useful in the proof of the following result.
\begin{lem}\label{homo}
    Let $M$ be a compact manifold without conjugate points and with visibility universal covering $\Tilde{M}$. Then for every $\theta\in T_1\Tilde{M}$ the projection map $P_{\theta}$ is a well-defined homeomorphism onto its image such that $I(\theta)\subset Im(P_{\theta})$.    
\end{lem}

\begin{proof}
As above we denote by $\eta$ the unique unit vector induced by angles $(\alpha_1,\ldots ,\alpha_{n-1}) \in [0,a]^{n-1}$ and by $P_{\theta}(\eta)$ the intersection $\gamma_{\eta}\cap H^+(\theta)$. If we choose $\epsilon>0$ small enough in Theorem \ref{T} then $P_{\theta}(\eta)$ always exist. The intersection $\gamma_{\eta}\cap H^+(\theta)$ is unique because of the equidistance of horospheres (Proposition \ref{horo}(4)) so $P_{\theta}$ is well-defined. Since $\Tilde{M}$ has no conjugate points it follows that $P_{\theta}$ is injective. Let $\eta$ be induced by any angles in $Dom(P_{\theta})$ and $\epsilon_1>0$. Since $\exp_{\gamma_{\theta}(T)}$ is a diffeomorphism, for $\epsilon_1$ there exists $\delta>0$ such that if $\xi$ is $\delta$-close to $\eta$ then $d(\gamma_{\eta}(T),\gamma_{\xi}(T))\leq \epsilon_1$. Furthermore, Theorem \ref{T} provides that $d(\gamma_{\eta}(T),\gamma_{\eta}\cap H^+(\theta))\leq \epsilon$ and $d(\gamma_{\xi}(T),\gamma_{\xi}\cap H^+(\theta))\leq \epsilon$ for every $\xi,\eta\in V(\gamma_{\theta}(T))$ induced by any angles in $Dom(P_{\theta})$. Applying the triangle inequality we get the continuity of $P_{\theta}$: for every $\xi$ $\delta$-close to $\eta$,
$$d(P_{\theta}(\eta),P_{\theta}(\xi))\leq d(P_{\theta}(\eta),\gamma_{\eta}(T))+d(\gamma_{\eta}(T),\gamma_{\xi}(T))+d(\gamma_{\xi}(T),P_{\theta}(\xi))\leq 2\epsilon+\epsilon_1.$$
Thus $P_{\theta}$ is a continuous bijection from a compact space $[0,a]^{n-1}$ onto its Hausdorff image hence a homeomorphism. Since $Diam(I(\theta))\leq Q$ we conclude that $I(\theta)\subset H^+(\theta)\cap \bar{B}_{2Q}(p)$. Thus choosing $\epsilon>0$ small enough, Theorem \ref{T} ensures that $I(\theta)\subset Im(P_{\theta})$.
\end{proof}

\begin{cor}\label{par}
Let $M$ be a compact manifold without conjugate points and with visibility universal covering $\Tilde{M}$. Then for every $\theta\in T_1\Tilde{M}$, there exist continuous parametrizations of some open sets $V\subset H^+(\theta)$, $U\subset \mathcal{\Tilde{F}}^s(\theta)$ with $I(\theta)\subset V$ and $\mathcal{\Tilde{I}}(\theta)\subset U$ by $(n-1)$ polar coordinates. An analogous statement holds for the unstable case, i.e., $H^-(\theta)$ and $\mathcal{\Tilde{F}}^u(\theta)$.
\end{cor}
\begin{proof}
For every $\theta \in T_1\Tilde{M}$, Lemma \ref{homo} gives a homeomorphism between $[0,a]^{n-1}$ and $Im(P_{\theta})\subset H^+(\theta)$ containing $I(\theta)$. Choosing an open set $V\subset Im(P_{\theta})$ containing $I(\theta)$, we get a continuous parametrization of $V\subset H^+(\theta)$ by polar coordinates. Now, we attach to each $p\in V$ a unit vector $v_p\in T_p\Tilde{M}$ orthogonal to $H^+(\theta)$ and directed inward. Clearly, the set $U=\{ (p,v_p)\in T_1\Tilde{M}: p\in V \}$ is an open set of $\mathcal{\Tilde{F}}^s(\theta)$ containing $\mathcal{\Tilde{I}}(\theta)$. Since $V$ is parametrized by polar coordinates, so is $U$.
\end{proof}

\subsection{Some dynamical properties and measures of maximal entropy for continuous flows}
In this subsection we assume a continuous flow $\psi_t:X\to X$ acting on a compact metric space $X$. We state some dynamical and ergodic properties we will work with in later sections. Let us begin with the stable and unstable sets. For every $x\in X$, we define the strong stable set of $x$ by
$$W^{ss}(x)=\{ y\in X:d(\psi_t(x),\psi_t(y))\to 0 \text{ as }t\to \infty\},$$
and for every $\epsilon>0$, the $\epsilon$-strong stable set of $x$ by 
$$W^{ss}_{\epsilon}(x)=\{ y\in W^{ss}(x):d(\psi_t(x),\psi_t(y))\leq \epsilon \text{ for every }t\geq 0\}.$$
The strong unstable set $W^{uu}(x)$ and the $\epsilon$-strong unstable set $W^{uu}_{\epsilon}(x)$ are defined analogously for $t\leq0$. All the above sets are non-empty for Anosov and expansive homeomorphisms on surfaces and for Anosov and expansive flows on $3$-manifolds. We say that the flow $\psi_t$ is expansive if there exists $\epsilon>0$ such that if $x,y\in X$ satisfy $d(\psi_t(x),\psi_{\rho(t)}(y))\leq \epsilon$ for every $t\in \mathbb{R}$ and some continuous surjection $\rho:\mathbb{R} \to \mathbb{R}$ with $\rho(0)=0$, then there exists $\tau\in [-\epsilon,\epsilon]$ with $y=\psi_{\tau}(x)$. We call $\epsilon>0$ a constant of expansivity for $\psi_t$. In the context of continuous flows without singularities acting on compact manifolds, the above definition is equivalent to the Bowen-Walters definition \cite{bowen72}. 

The $\epsilon$-strong stable and $\epsilon$-strong unstable sets help in the definition of the so-called local product for flows. The flow $\psi_t$ has a local product structure if for every $\epsilon>0$ there exists $\delta>0$ such that if $x,y\in X$ satisfy $d(x,y)\leq \delta$ then there exists a unique $\tau\in \mathbb{R}$ with 
$$|\tau|\leq \epsilon \quad \text{ and } \quad W^{ss}_{\epsilon}(x)\cap W^{uu}_{\epsilon}(\psi_{\tau}(y))\neq \emptyset.$$
We observe that the intersection is not unique in general and the intersection points accompany $x$ in the future and $y$ in the past. Anosov flows are typical examples which have local product structure. Furthermore, expansive homeomorphisms on surfaces have connected stable and unstable sets with local product structure everywhere except at a finite set of periodic orbits \cite{hir89,lew89}.

Finally, let us define a special kind of semi-conjugacy between flows. Let $\phi_t:Y\to Y$ and $\psi_t:X\to X$ be two continuous flows acting on compact topological spaces. A map $\chi:Y\to X$ is called a time-preserving semi-conjugacy if $\chi$ is a continuous surjection satisfying $\chi\circ\phi_t=\psi_t\circ \chi$ for every $t\in \mathbb{R}$. In this case, we say that $\psi_t$ is time-preserving semi-conjugate to $\phi_t$ or is a time-preserving factor of $\phi_t$.

We now look at measures of maximal entropy. A Borel set $Z\subset X$ is invariant by the flow if $\psi_t(Z)=Z$ for every $t\in \mathbb{R}$. A probability measure $\nu$ on $X$ is invariant by the flow if $(\psi_t)_*\nu=\nu$ for every $t\in \mathbb{R}$. Denote by $\mathcal{M}(\psi)$ the set of all flow-invariant-measures on $X$. A measure $\nu\in \mathcal{M}(\psi)$ is ergodic if for every flow-invariant set $A\subset X$, we have either $\nu(Z)=0$ or $\nu(Z)=1$.

Let $Z\subset X$ be a flow-invariant Borel set and $\nu$ be a flow-invariant measure supported on $Z$. We define the metric entropy $h_{\nu}(\psi,Z)$ of $\nu$ with respect to the flow $\psi$ as the metric entropy $h_{\nu}(\psi_1,Z)$ with respect to its time-1 map $\psi_1$ \cite{walt00}. For $Z=X$ we write $h_{\nu}(\psi)$. When $Z$ is also compact, we define the topological entropy of $Z$ as follows. For every $\epsilon,T>0$ and every $x\in Z$, we define the $(T,\epsilon)$-dynamical balls by
\[  B(x,\epsilon,T)=\{ y\in Z:d(\psi_s(x),\psi_s(y))<\epsilon,s\in [0,T] \}  \] 
Denote by $M(T,\epsilon,Z)$ the minimum cardinality of any cover of $Z$ by $(T,\epsilon)$-dynamical balls. The topological entropy of $Z$ with respect to $\psi$ is
\[  h(\psi,Z)=\lim_{\epsilon\to 0}\limsup_{T\to \infty} \frac{1}{T} \log M(T,\epsilon,Z).\]
For $Z=X$ we write $h(\psi)$. We remark that $h(\psi,Z)=h(\psi_1,Z)$ where $h(\psi_1,Z)$ is the topological entropy of $Z$ with respect to the time-1 map $\psi_1$. Observe that expansive flows are examples of continuous systems having positive topological entropy.

A fundamental result that relates the topological and metric entropy of each compact set $Z\subset X$ invariant by the flow $\psi_t$ is the so-called the variational principle \cite{dina71}:
\begin{equation}\label{varcon}
h(\psi,Z)= \sup_{\nu} h_{\nu}(\psi,Z),    
\end{equation}
where $\nu$ varies over all flow-invariant measures supported on $Z$. We say that $\mu\in \mathcal{M}(\psi)$ supported on $Z$ is a measure of maximal entropy if $h_{\mu}(\psi,Z)$ achieves the maximum in \eqref{varcon}. If $Z=X$ and $\mu$ is the only measure satisfying this condition then $\mu$ is called the unique measure of maximal entropy for the flow $\psi$. We will talk more about the role of this measure in our work in a later section.

\section{The factor flow}\label{c26}
In this section we extend the construction of the factor flow of the geodesic flow introduced by Gelfert and Ruggiero for the case of compact higher genus surfaces without focal points \cite{gelf19}. We always consider a compact $n$-manifold $M$ without conjugate points and with visibility universal covering $\tilde{M}$. 

Let us observe that in our setting, intersections $\mathcal{\Tilde{I}}$ still have good topological behaviour, i.e., intersections $\mathcal{I}(\theta)$ are compact connected sets uniformly bounded by Morse's constant $Q>0$ by Corollary \ref{morse} which is a consequence of the generalized flat strip Theorem.

We outline the construction of the factor flow of Section 4 of \cite{gelf19} in our context. For every $\eta,\theta\in T_1M$, $\eta$ and $\theta$ are equivalent, 
$$\eta\sim\theta \quad \text{ if and only if }\quad  \eta\in \mathcal{I}(\theta).$$
This is an equivalence relation that induces a quotient space $X$ and a quotient map $\chi:T_1M\to X$. For every $\theta \in T_1M$, we denote by $[\theta]=\chi(\theta)$ the equivalence class of $\theta$. Using the geodesic flow $\phi_t$ induced by $(M,g)$, we define a quotient flow $\psi_t:X\to X$ by the formula $\psi_t[\theta]=[\phi_t(\theta)]$ valid for every $t\in \mathbb{R}$. We shall endow $X$ with the quotient topology.

Similarly, we can repeat the above construction in the context of the universal covering $\tilde{M}$. We say that $\eta,\theta\in T_1\tilde{M}$ are equivalent if and only if $\eta \in\mathcal{\tilde{I}}(\theta)=\mathcal{\tilde{F}}^s(\theta)\cap \mathcal{\tilde{F}}^u(\theta)\subset T_1\tilde{M}$. This equivalence relation induces a quotient space $\tilde{X}$, a quotient map $\Tilde{\chi}:T_1\tilde{M}\rightarrow \tilde{X}$ and a quotient flow $\tilde{\psi}_t:\tilde{X}\rightarrow \Tilde{X}$. Since $T_1\tilde{M}$ is the universal covering of $T_1M$, we can show the map $\Pi:\tilde{X}\to X$ satisfying $\chi\circ d\pi =\Pi\circ \tilde{\chi}$ is a well-defined covering map. Thus, $\tilde{X}$ and $X$ are locally homeomorphic when both are endowed with the quotient topology.

Before stating some properties of the quotient space and quotient flow we introduce some useful concepts and results. For every subset $A\subset T_1M$, the \textbf{saturation} of $A$ with respect to the quotient map $\chi$ is defined by $Sat(A)=\chi^{-1}\circ \chi(A)$. This can also be expressed by
$$Sat(A)=\bigcup_{\eta\in A}\chi^{-1}\circ \chi(\eta)=\bigcup_{\eta\in A}Sat(\eta).$$ 
From this formula, since $\eta\in Sat(\eta)$, we conclude that any $Sat(\eta)\subset Sat(A)$ cannot be disjoint of $A$. A subset $A\subset T_1M$ is \textbf{saturated} with respect to $\chi$ if $A=Sat(A)$. By the above formula, a saturated set contains all the saturations of its elements. Furthermore, $Sat(A)$ is the smallest saturated set containing $A$.

Taking into account the above equivalence relation, we see that for every $\eta\in T_1M$ and every $A\subset T_1M$,
$$Sat(\eta)=\chi^{-1}\circ \chi(\eta)=\mathcal{I}(\eta) \quad \text{ and }\quad Sat(A)=\bigcup_{\eta\in A}\mathcal{I}(\eta).$$
So, every expansive point $\eta\in T_1M$ has $Sat(\eta)=\eta$ and every non-expansive point does not have a singleton saturated set. We say that $\eta$ has a trivial equivalence class if $Sat(\eta)=\eta$ otherwise $\eta$ has a non-trivial equivalence class. The following result provides a way to construct open saturated sets in $T_1M$.
\begin{lem}\label{sat}
Let $M$ be a compact manifold without conjugate points and with visibility universal covering $\tilde{M}$. Then, 
\begin{enumerate}
    \item If $C\subset T_1M$ is a closed set then so is $Sat(C)$.
    \item If $U\subset T_1M$ is an open set then $U'=T_1M\setminus Sat(T_1M\setminus U)$ is an open saturated set contained in $U$. In particular, if $U$ contains some $\mathcal{I}(\eta)$, so does $U'$.
\end{enumerate}
\end{lem}
\begin{proof}
For item 1, let $\theta\in T_1M$ be an accumulation point of $Sat(C)$ and $(\theta_n)\subset Sat(C)$ be a sequence converging to $\theta$. We have to show that $\theta \in Sat(C)$. Without loss of generality we can suppose that $(\theta_n)$ has a subsequence $(\theta_{n_k})$ composed only of non-expansive points. We claim that for every $k\geq 1$ there exists $\xi_{n_k}\in Sat(\theta_{n_k})\cap C$. Otherwise for some $k\geq 1$, $Sat(\theta_{n_k})$ is disjoint of $C$, a contradiction because $Sat(\theta_{n_k})\subset Sat(C)$. Since $\theta_{n_k}\to \theta$ and $\xi_{n_k}\in Sat(\theta_{n_k})=\mathcal{I}(\theta_{n_k})=\mathcal{F}^s(\theta_{n_k})\cap \mathcal{F}^u(\theta_{n_k})$, the continuity of the horospherical foliations $\mathcal{F}^s$ and $\mathcal{F}^u$ implies that $\xi_{n_k}$ converges to some $\xi\in \mathcal{F}^s(\theta)\cap \mathcal{F}^u(\theta) =\mathcal{I}(\theta)$. Recalling that $(\xi_{n_k})\subset C$ it follows that $\xi\in C$ hence $\theta\in \mathcal{I}(\theta)=\mathcal{I}(\xi)=Sat(\xi)\subset Sat(C)$.
In item 2, applying item 1 to $T_1M\setminus U$ we get that $Sat(T_1M\setminus U)$ is closed hence $U'$ is an open set. It is clear to verify that $U'\subset U$. Observing that complements of saturated sets are also saturated, we obtain the result.
\end{proof}
Open saturated sets in $T_1M$ are quite useful because their images under $\chi$ are open sets in $X$ endowed with the quotient topology. We apply this fact to give a short proof of the following property.
\begin{lem}\label{haus}
The quotient space $X$ is Hausdorff.
\end{lem}
\begin{proof}
Let $[\xi],[\eta]\in X$ be distinct hence $\mathcal{I}(\xi)$ and $\mathcal{I}(\eta)$ are disjoint compact sets. Thus, there exist disjoint open neighborhoods $U,V\subset T_1M$ of $\mathcal{I}(\xi)$ and $\mathcal{I}(\eta)$. By Lemma \ref{sat}(2), $U'=T_1M\setminus Sat(T_1M\setminus U)\subset U$ and $V'=T_1M\setminus Sat(T_1M\setminus V)\subset V$ are disjoint open saturated neighborhoods of $\mathcal{I}(\xi)$ and $\mathcal{I}(\eta)$. So, $\chi(U')$ and $\chi(V')$ are disjoint open neighborhoods of $X$ separating $[\xi]$ and $[\eta]$.
\end{proof}
Another tool we will use is the following basic proposition of general topology.
\begin{pro}[\cite{will12}]\label{topol}
If $f:X\to Y$ is a continuous surjection from a compact metric space onto a Hausdorff space then $Y$ is metrizable.
\end{pro}
We now state the basic properties of the quotient space and quotient flow.
\begin{lem}\label{factor}
Let $M$ be a compact $n$-dimensional manifold without conjugate points and with visibility universal covering $\Tilde{M}$, $\phi_t$ be its geodesic flow, $\psi_t$ be the quotient flow acting on the quotient space $X$ and $\chi$ be the quotient map. Then,
\begin{enumerate}
    \item $\psi_t$ is a continuous flow time-preserving semi-conjugate to $\phi_t$ under the quotient map $\chi$. In particular, $\psi_t$ is a factor flow of the geodesic flow.
    \item $X$ is a compact metric space of topological dimension at least $n-1$.
\end{enumerate}
\end{lem}
\begin{proof}
For item 1, $\psi_t$ is well-defined because the horospherical foliations are invariant by the geodesic flow. Using the definition formula of $\psi_t$ and the quotient topology it is clear to verify that $\psi_t$ is a continuous flow acting on $X$. Moreover, $\psi_t$ is time-preserving semi-conjugate to $\phi_t$ by definition.

For item 2, $X$ is compact since $\chi$ is continuous. Applying Proposition \ref{topol} and Lemma \ref{haus} to the quotient map $\chi:T_1M\to X$ we get that $X$ is metrizable. For the topological dimension the argument in \cite{mam23} extends naturally. Let $\theta \in T_1\tilde{M}$ and $V(\theta)$ be its vertical fiber. We define the set $W_{\theta}=\bigcup_{t\in \mathbb{R}}\tilde{\phi}_t(V(\theta))\subset T_1\tilde{M}$. Since $V(\theta)$ is homeomorphic to the sphere $\mathbb{S}^{n-1}$, we see that $W_{\theta}$ is a topological hypersurface of dimension $n-1$. By the divergence of geodesic rays, for every $\eta,\xi\in V(\theta)$ it holds that $\eta\not\in \mathcal{\tilde{I}}(\xi)$. So, the restriction of $\tilde{\chi}$ to $W_{\theta}$ is injective hence bijective onto its image. This implies that $\tilde{\chi}(W_{\theta})\subset \tilde{X}$ is homeomorphic to a hypersurface of dimension $n-1$ and so the topological dimension of $\Tilde{X}$ is at least $n-1$. This conclusion extends to $X$ because $\tilde{X}$ and $X$ are locally homeomorphic.
\end{proof}
Using this lemma we can define a distance on $X$, denoted by $d$, which is compatible with the quotient topology. We highlight that except for the compactness and metrizability of $X$, analogous definitions and properties hold for objects in the setting of the universal covering $\Tilde{M}$ such as $\Tilde{X}, \Tilde{\chi}, \tilde{X}$ and the open sets of $T_1\Tilde{M}$.

\section{Basis of neighborhoods of the quotient space}\label{basis}
The goal of this section is to construct a special basis of neighborhoods for the quotient topology. In \cite{gelf19}, Gelfert and Ruggiero constructed a similar basis of neighborhoods for the surface case. This construction cannot be extended to our setting due to the topological and geometrical difficulties in higher dimensions. We will make another type of construction using the parametrizations of the horospherical leaves (Section \ref{parame}). We will assume a compact manifold $M$ without conjugate points and with visibility universal covering $\tilde{M}$. 

Choose an arbitrary $[\theta]\in \Tilde{X}$ for some $\theta\in T_1\Tilde{M}$. Applying Corollary \ref{par} to $\theta$ we get an open set $U\subset \mathcal{\Tilde{F}}^s(\theta)$ parametrized by $(n-1)$ polar coordinates satisfying $\mathcal{\Tilde{I}}(\theta)\subset U$. So, $U$ is homeomorphic to an open set of $\mathbb{R}^{n-1}$ and hence we can pullback the distance of $\mathbb{R}^{n-1}$ to $U$ and denote it by $d_U$. Since $\mathcal{\Tilde{I}}(\theta)\subset U$, for each $\delta>0$ small enough consider the $\delta$-open neighborhood of $\mathcal{\Tilde{I}}(\theta)$
$$W_{\delta}(\mathcal{\Tilde{I}}(\theta))=\{ \xi\in U: d_U(\xi,\eta)<\delta \text{ for some } \eta \in \mathcal{\Tilde{I}}(\theta) \}.$$
Now, for every $\epsilon>0$ and every $\eta\in T_1\Tilde{M}$ denote by $V_{\epsilon}(\eta)$ the open $\epsilon$-ball of the vertical subspace $V(\eta)$, which is centered at $\eta$. For $t'>0$ small enough we define the following set 
$$U^s_{\epsilon,\delta,t'}(\theta)=\bigcup_{|t|<t', \eta\in W_{\delta}(\mathcal{\Tilde{I}}(\theta))} \Tilde{\phi}_t(V_{\epsilon}(\eta)).$$
Since both sets $W_{\delta}(\mathcal{\Tilde{I}}(\theta))$ and $V_{\epsilon}(\eta)$ are parametrized by $(n-1)$ variables each, we see that $U^s_{\epsilon, \delta, t'}(\theta)$ can be continuously parametrized by $(2n-1)$ variables and hence $U^s_{\epsilon, \delta, t'}(\theta)$ is the image of an open set of $\mathbb{R}^{2n-1}$ by a continuous injective map. Applying the invariant domain Theorem we get that $U^s_{\epsilon, \delta, t'}(\theta)$ is an open set of $T_1\tilde{M}$ containing $\mathcal{\Tilde{I}}(\theta)$.

In this setting, we define an open section transversal to the geodesic flow
$$\Sigma^s_{\epsilon,\delta}(\theta)=\bigcup_{\eta\in V_{\epsilon}(\theta)} (\mathcal{\Tilde{F}}^s(\eta)\cap U^s_{\epsilon, \delta, t'}(\theta)).$$
Indeed $\Sigma^s_{\epsilon,\delta}(\theta)$ is a set composed by pieces of stable horospherical leaves $\mathcal{\tilde{F}}^s(\eta)$ which are orthogonal to the geodesic flow. Therefore $\Sigma^s_{\epsilon,\delta}(\theta)$ is a section transversal to the geodesic flow. Since $\eta$ varies continuously over $V_{\epsilon}(\theta)$, the continuity of the horospherical foliations implies that $\Sigma^s_{\epsilon,\delta}(\theta)$ is an open transversal section. Furthermore, it holds that $\mathcal{\Tilde{I}}(\theta)\subset \Sigma^s_{\epsilon,\delta}(\theta)$. Choosing $\tau>0$ small enough we write
$$U^s_{\epsilon, \delta,\tau}(\theta)= \bigcup_{|t|<\tau} \Tilde{\phi}_t(\Sigma^s_{\epsilon,\delta}(\theta)).$$
Clearly, $U^s_{\epsilon, \delta,\tau}(\theta)$ is an open set of $T_1\Tilde{M}$ containing $\mathcal{\Tilde{I}}(\theta)$. Moreover, $U^s_{\epsilon, \delta,\tau}(\theta)$ is foliated by pieces of center stable sets $\mathcal{\Tilde{F}}^{cs}(\eta)$ for $\eta\in \Sigma^s_{\epsilon,\delta}(\theta)$. Observe that diameter of $U^s_{\epsilon, \delta,\tau}(\theta)$ is controlled by $\epsilon$ in the direction of $V_{\epsilon}(\theta)$, by $\delta$ in the direction of $\mathcal{\Tilde{F}}^s(\theta)$ and by $\tau$ in the geodesic flow direction. Thus, when $\epsilon,\delta,\tau\to 0$, diameter of $U^s_{\epsilon, \delta,\tau}(\theta)$ tends to diameter of $\mathcal{\Tilde{I}}(\theta)$. 

Since Corollary \ref{par} also holds for the unstable case, we can make an analogous construction using unstable horospherical leaves and the same parameters $\epsilon,\delta,\tau,t'>0$. We thus get another open set $U^u_{\epsilon, \delta,\tau}(\theta)\subset T_1\Tilde{M}$ containing $\mathcal{\Tilde{I}}(\theta)$, using analogous formulas 
$$\Sigma^u_{\epsilon,\delta}(\theta)=\bigcup_{\eta\in V_{\epsilon}(\theta)} (\mathcal{\Tilde{F}}^u(\eta)\cap U^u_{\epsilon, \delta, t'}(\theta)) \quad\text{ and }\quad U^u_{\epsilon, \delta,\tau}(\theta)= \bigcup_{|t|<\tau} \Tilde{\phi}_t(\Sigma^u_{\epsilon,\delta}(\theta)).$$
This time $U^u_{\epsilon, \delta,\tau}(\theta)$ is foliated by pieces of center unstable sets $\mathcal{F}^{cu}(\eta)$ for $\eta \in \Sigma^u_{\epsilon,\delta}(\theta)$. Again, the diameter of $U^u_{\epsilon, \delta,\tau}(\theta)$ is controlled by $\epsilon$ in the direction of $V_{\epsilon}(\theta)$, by $\delta$ in the direction of $\mathcal{\Tilde{F}}^u(\theta)$ and by $\tau$ in the geodesic flow direction. 

So far, we have two open neighborhoods of $\mathcal{\Tilde{I}}(\theta)$, $U^s_{\epsilon, \delta,\tau}(\theta)$ foliated by pieces of center stable sets and $U^u_{\epsilon, \delta,\tau}(\theta)$ foliated by pieces of center unstable sets. Using the same parameter values $\epsilon,\delta,\tau>0$, we take the intersection of these open neighborhoods
$$U_{\epsilon, \delta,\tau}(\theta)=U^s_{\epsilon, \delta,\tau}(\theta)\cap U^u_{\epsilon, \delta,\tau}(\theta).$$
Clearly, $U_{\epsilon, \delta,\tau}(\theta)$ is an open set of $T_1\Tilde{M}$ containing $\mathcal{\Tilde{I}}(\theta)$. Moreover, $U_{\epsilon, \delta,\tau}(\theta)$ is composed of intersections of stable and unstable horospherical leaves. By above, the diameter of $U_{\epsilon, \delta,\tau}(\theta)$ is still controlled by parameters $\epsilon,\delta,\tau>0$.

\begin{figure}[ht]
    \centering
    \includegraphics[width=\textwidth]{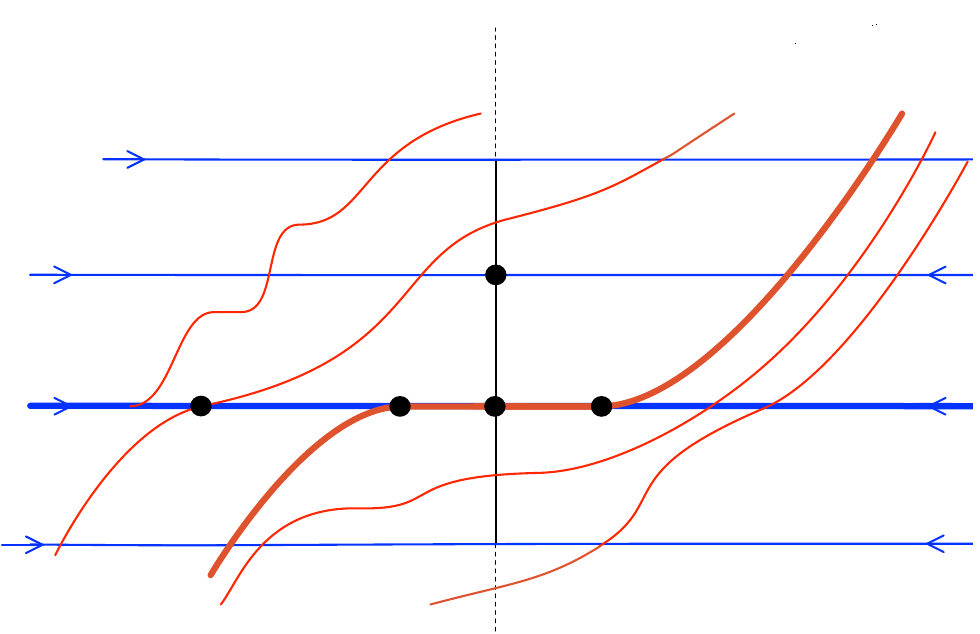}
    \caption{Some cross section of $A_{\epsilon, \delta,\tau}(\theta)$ and its properties for the surface case (figure borrowed from K. Gelfert \cite{gelf20}).}
    \label{kato}
\end{figure}

Let us observe that $U_{\epsilon, \delta,\tau}(\theta)$ is not saturated in general. So, applying Lemma \ref{sat} to $U_{\epsilon, \delta,\tau}(\theta)$ gives rise to a smaller open saturated set $A_{\epsilon, \delta,\tau}(\theta)\subset U_{\epsilon, \delta,\tau}(\theta)$. Since $\mathcal{\Tilde{I}}(\theta)\subset U_{\epsilon, \delta,\tau}(\theta)$, by definition we see that $\mathcal{\Tilde{I}}(\theta)\subset A_{\epsilon, \delta,\tau}(\theta)$ hence $Diam(A_{\epsilon, \delta,\tau}(\theta))$ is still controlled by parameters $\epsilon,\delta,\tau>0$. We thus obtain an open saturated set $A_{\epsilon, \delta,\tau}(\theta)$ containing $\mathcal{\Tilde{I}}(\theta)$, composed of intersections of stable and unstable horospherical leaves. Moreover, every two points of $A_{\epsilon, \delta,\tau}(\theta)$ are heteroclinically related by Theorem \ref{conec}, because $\Tilde{M}$ is a visibility manifold (see figure \ref{kato}). Varying $\epsilon,\delta,\tau$ we get a family of open neighborhoods of $\mathcal{\Tilde{I}}(\theta)$ that provides a basis of neighborhoods of $[\theta]\in \Tilde{X}$ in the quotient topology.

\begin{lem}\label{que1}
For every $\theta\in T_1\tilde{M}$, the family
$$\mathcal{A}_{\theta}= \{ \tilde{\chi}( A_{\epsilon_l, \delta_m,\tau_n}(\theta)): \epsilon_l=1/l,\delta_m=1/m,\tau_n=1/n \text{ with } l,m,n\in \mathbb{N} \}$$
is a countable basis of neighborhoods of $[\theta]\in \tilde{X}$. Hence $\tilde{X}$ is first countable and $\{ \mathcal{A}_{\theta}: \theta\in T_1\tilde{M} \}$ is a basis for the quotient topology of $\tilde{X}$. 
\end{lem}
\begin{proof}
For every $\theta\in T_1\tilde{M}$, $A=A_{\epsilon, \delta,\tau}(\theta)$ is a saturated open neighborhood of $\mathcal{\tilde{I}}(\theta)$ hence $\tilde{\chi}^{-1}\circ \tilde{\chi}(A)=A$. Therefore $\tilde{\chi}(A)$ is an open set of $\tilde{X}$ containing $[\theta]$. Choosing $\epsilon_l=1/l,\delta_m=1/m,\tau_n=1/n$ with $l,m,n\in \mathbb{N}$ large enough, gives a countable family $\mathcal{A}_{\theta}=\{ \tilde{\chi}( A_{\epsilon_l, \delta_m,\tau_n}(\theta))\}$ of open neighborhoods of $[\theta]\in \Tilde{X}$. 

For the basis property, let us first show that $\bigcap_{l,m,n\in \mathbb{N}}A_{\epsilon_l, \delta_m,\tau_n}(\theta)=\mathcal{\Tilde{I}}(\theta)$. The reverse inclusion is clear, for the direct one let $\eta\in \bigcap_{l,m,n\in \mathbb{N}}A_{\epsilon_l, \delta_m,\tau_n}(\theta)$. Recall that when $l,m,n\to \infty$, $Diam(A_{\epsilon_l, \delta_m,\tau_n}(\theta))\to Diam(\mathcal{\Tilde{I}}(\theta))$ and hence $Diam(A_{\epsilon_l, \delta_m,\tau_n}(\theta)\setminus \mathcal{\Tilde{I}}(\theta))\to 0$. Since $\eta\in A_{\epsilon_l, \delta_m,\tau_n}(\theta)$ for every $l,m,n\in \mathbb{N}$ large enough, the continuity of Sasaki distance implies that $d_s(\eta,\mathcal{\Tilde{I}}(\theta))\leq\lim_{l,m,n\to\infty}Diam(A_{\epsilon_l, \delta_m,\tau_n}(\theta)\setminus \mathcal{\Tilde{I}}(\theta))=0$. Since $\mathcal{\Tilde{I}}(\theta)$ is closed we conclude that $\eta \in\mathcal{\Tilde{I}}(\theta)$ which proves the direct inclusion. Now, applying the quotient map $\Tilde{\chi}$ to $\bigcap_{l,m,n\in \mathbb{N}}A_{\epsilon_l, \delta_m,\tau_n}(\theta)=\mathcal{\Tilde{I}}(\theta)$ yields $\bigcap_{l,m,n\in \mathbb{N}}\tilde{\chi}(A_{\epsilon_l, \delta_m,\tau_n}(\theta))=[\theta]$. This implies that $\mathcal{A}_{\theta}$ is basis of open neighborhoods for $[\theta]\in \Tilde{X}$.  
\end{proof} 
So far, we have a countable family of open neighborhoods for every $\mathcal{\tilde{I}}(\theta)\subset T_1\tilde{M}$ and a countable basis of open neighborhoods for every $[\theta]\in \tilde{X}$. Projecting by the covering maps $\Pi$ and $d\pi$, we get corresponding families of open neighborhoods of every $\mathcal{I}(\theta)\subset T_1M$ and $[\theta]\in X$. Therefore, $X$ is first countable and $\{ \Pi(\mathcal{A}_{\theta}): \theta\in T_1M \}$ is a basis for the quotient topology of $X$. 

\section{Dynamical properties of the factor flow}\label{dyn}
This section deals with some dynamical properties of the factor flow. We mention that properties follow an analogous development of \cite{mam23}.

Let us first state the theorem with the properties we are interested in.
\begin{teo}\label{prop}
Let $M$ be a compact manifold without conjugate points and with visibility universal covering $\Tilde{M}$, $\psi_t$ be the factor flow. Then, $\psi_t$ is topologically mixing, expansive and has a local product structure and the specification property. 
\end{teo}
The proof follows verbatim arguments to the case of compact higher genus surfaces without conjugates points, i.e., Section 6 in \cite{mam23}. So, we omit the proof and only restate some important facts. A key result is the following lemma because it allows us to focus only in the metric properties of $\psi_t$ (given by the distance) instead of the dimension of $M$. Recall that $d_s$ denotes the Sasaki distance which is defined in Equation \eqref{sasa}. 
\begin{lem}\label{lev}
Let $\Tilde{Q}>0$ be the Morse's constant given in Corollary \ref{morse}. Then, there exist $r_0,s_0>0$ such that for every $[\xi],[\eta]\in X$ with $d([\xi],[\eta])\leq r_0$ then
$$d_s(\tilde{\xi},\tilde{\eta})\leq \Tilde{Q}+s_0,$$ 
for some lifts $\tilde{\xi},\tilde{\eta}\in T_1\tilde{M}$ of $\xi,\eta\in T_1M$.
\end{lem}
\begin{proof}
Denote by $\Tilde{A}_{\epsilon,\delta,\tau}(\Tilde{\theta})$ the open neighborhood of $\mathcal{\Tilde{I}}(\Tilde{\theta})$ given in Lemma \ref{que1} and by $A_{\epsilon,\delta,\tau}(\theta)$ the projection of $\Tilde{A}_{\epsilon,\delta,\tau}(\Tilde{\theta})$ by the covering map $d\pi:T_1\Tilde{M}\to T_1M$. For every $\theta\in T_1M$, choose $\epsilon,\delta,\tau>0$ small enough so that $A_{\epsilon,\delta,\tau}(\theta)$ is evenly covered by $d\pi$. Since the family $\mathcal{A}=\{ \chi(A_{\epsilon,\delta,\tau}(\theta)): \theta \in T_1M\}$ is an open covering of the compact space $X$, there exists a Lebesgue number $r_0>0$ for $\mathcal{A}$. So, for every $[\eta],[\xi]\in X$ with $d([\eta],[\xi])\leq r_0$, there exists $\theta\in T_1M$ such that $[\xi]\in B([\eta],r_0)\subset \chi(A_{\epsilon,\delta,\tau}(\theta))\in \mathcal{A}$ where $B([\eta],r_0)$ is the closed ball of radius $r_0$ centered at $[\eta]$. Since $A_{\epsilon,\delta,\tau}(\theta)$ is evenly covered by $d\pi$, it follows that for any lift $\tilde{\theta}\in T_1\Tilde{M}$ of $\theta$ there exist lifts $\tilde{A}_{\epsilon,\delta,\tau}(\tilde{\theta}),\tilde{\eta},\tilde{\xi}$
of $A_{\epsilon,\delta,\tau}(\theta),\eta,\xi$ such that $\tilde{\eta},\tilde{\xi}\in \tilde{A}_{\epsilon,\delta,\tau}(\tilde{\theta})\subset T_1\Tilde{M}$. By our choice of $\epsilon,\delta,\tau$ we know that there exists $s_0>0$ such that $Diam(\tilde{A}_{\epsilon,\delta,\tau}(\tilde{\theta}))\leq Q+s_0$ and therefore $d_s(\tilde{\eta},\tilde{\xi})\leq Q+s_0$.
\end{proof}
The proof of expansivity of the factor flow is the same as in the surface case. We give it again below because it is one of the most important results of the work. Its importance lies in the fact that all subsequent results essentially rely on this property.
\begin{proof}
Let $r_0>0$ be given by Lemma \ref{lev}. Given two orbits of $\psi_t$ having Hausdorff distance bounded by $r_0$, we will show that the orbits agree. Let $[\eta],[\xi]\in X$ with $d(\psi_t[\eta],\psi_{\rho(t)}[\xi])\leq r_0$ for every $t\in \mathbb{R}$ and some reparametrization $\rho$. By Lemma \ref{lev}, there exist lifts $\tilde{\eta},\tilde{\xi}\in T_1\tilde{M}$ of $\eta,\xi$ such that $d_s(\phi_t(\tilde{\eta}),\phi_{\rho(t)}(\tilde{\xi}))\leq Q+s_0$ for every $t\in \mathbb{R}$. Thus, the orbits of $\tilde{\eta}$ and $\tilde{\xi}$ have Hausdorff distance bounded by $Q+s_0$ hence the orbits are bi-asymptotic. This implies that there exists $\tau\in \mathbb{R}$ so that $\tilde{\xi}\in \mathcal{\tilde{I}}(\phi_{\tau}(\tilde{\eta}))$ hence $[\xi]=\psi_{\tau}[\eta]$. It is not hard to show that $|\tau|\leq \epsilon$ for some $\epsilon>0$ and so $\min(r_0,\epsilon)$ is an expansivity constant for $\psi_t$.
\end{proof}
For the local product structure, natural candidates for being the stable and unstable sets of $\psi_t$ are the images by $\chi$ of the stable and unstable horospherical leaves. For every $\eta\in T_1M$, we define the following sets
$$V^*[\eta]=\chi(\mathcal{F}^*(\eta)), \qquad *=s,u,cs,cu.$$
The heteroclinic connections (Theorem \ref{conec}) and the construction of open saturated sets (Section \ref{basis}) imply that any distinct $\eta,\xi\in A_{\epsilon_l, \delta_m,\tau_n}(\theta)$ are heteroclinically related hence so are $[\eta],[\xi]\in \chi(A_{\epsilon_l, \delta_m,\tau_n}(\theta))$: 
$$V^s[\eta]\cap V^{cu}[\xi]=V^s[\eta]\cap V^u(\psi_{\tau}[\xi])=\chi(\mathcal{F}^s(\eta)\cap \mathcal{F}^u(\phi_{\tau}(\xi))).$$
For $*=s,u,cs,cu$, $V^*[\eta]$ denotes the connected component containing $[\eta]$ and included in $\chi(A_{\epsilon_l, \delta_m,\tau_n}(\theta))$. The remaining lemmas in \cite{mam23} focus on proving the local product structure precisely, to which we refer for the proof. We join in the following proposition some important intermediate results.
\begin{pro}\label{loc}
Let $r_0>0$ be given by Lemma \ref{lev}. Then,
\begin{enumerate}
    \item For every $\epsilon>0$ and $D\in(0,r_0]$ there exists $T>0$ such that if $[\eta]\in V^s[\xi]\cap B([\xi],r_0)$ and $d([\eta],[\xi])\leq D$ then $d(\psi_t[\eta],\psi_t[\xi])\leq \epsilon$ for every $t\geq T$ (Uniform contraction).
    \item For every $[\eta]\in X$ and every $\epsilon>0$ there exists $\delta\in (0,r_0]$ such that $W^{ss}[\eta]\cap B([\eta],r_0)=V^s[\eta]\cap B([\eta],r_0)$ and $W^{ss}_{\epsilon}[\eta]\cap B([\eta],\delta)=V^s[\eta]\cap B([\eta],\delta)$.
\end{enumerate}
In particular, $V^s[\eta]$ agrees with the strong stable $W^{ss}[\eta]$ and $\epsilon$-strong stable $W_{\epsilon}^{ss}[\eta]$ sets locally. Analogous results hold for the unstable case.   
\end{pro}
For the specification property definition, we remit to Section 6 of \cite{gelf19}. In our setting this property is also obtained from a combination of some dynamical properties of $\psi_t$: topological mixing, expansivity and the local product structure. The specification property is also important because of the following classical result of Bowen and Franco \cite{bowen71endo,fran77}.
\begin{teo}\label{frank}
Let $\phi_t:X\to X$ be a continuous flow acting on a compact metric space. If $\phi_t$ is expansive and has the specification property then $\phi_t$ has a unique measure of maximal entropy.
\end{teo}
Applying this theorem to the factor flow $\psi_t$, we get the unique measure of maximal entropy $\nu$ for $\psi_t$. This measure will be useful for obtaining the measure of maximal entropy for the geodesic flow $\phi_t$.

\section{Uniqueness of the measure of maximal entropy}\label{ent}
In this section we extend the uniqueness of the measure of maximal entropy for the geodesic flow from the surface case to the higher dimensional case. To show the uniqueness we follow the same strategy as \cite{mam23}. As mentioned in Section \ref{dyn} the main difficulty was to extend some key lemmas to the higher dimensional setting, in particular, the construction of a special basis of neighborhoods of the quotient topology. 

For compact manifolds of negative curvature, Margulis \cite{margu70} constructed a measure using the stable and unstable invariant manifolds. Later, Bowen showed that this measure is actually a measure of maximal entropy for the geodesic flow. Using symbolic dynamics Bowen \cite{bowenper,bowen73maxi} also proved the uniqueness of this measure. In 1998, Knieper \cite{knie98} extended the existence and uniqueness to compact rank-1 manifolds of non-positive curvature. In 2016, Climenhaga and Thompson \cite{clim16} showed a non-uniform version of a classical theorem due to Bowen and Franco \cite{fran77}. This theorem helps to show the uniqueness of measures of maximal entropy for certain systems. In 2018, Gelfert and Ruggiero \cite{gelf19} proved the existence and uniqueness for compact higher genus surfaces without focal points using an expansive factor flow of the geodesic flow. In 2021, Climenhaga, Knieper and War \cite{clim21} extended this result to the higher dimensional setting which includes the case of compact higher genus surfaces without conjugate points. To show the uniqueness, they used Climenhaga-Thompsom non-uniform criterion. We restate Climehaga-Knieper-War's Theorem to compare it with our contribution. Recall that $\mathcal{R}_0$ is the set of expansive points.
\begin{teo}\label{clime}
    Let $(M,g)$ be a compact $n$-manifold without conjugate points and $\phi_t$ be its geodesic flow. Assume that
    \begin{enumerate}
        \item $M$ supports a metric $g_0$ of negative curvature.
        \item Geodesic rays diverge in the universal covering $\Tilde{M}$.
        \item The fundamental group $\pi_1(M)$ is residually finite.
        \item $\sup \{ h_{\mu}(\phi_1): \mu \text{ is supported on } T_1M\setminus \mathcal{R}_0 \}<h(\phi_1)$.
    \end{enumerate}
    Then, the geodesic flow has a unique measure of maximal entropy which has full support.
\end{teo}
In particular, this conclusion holds for compact higher genus surfaces without conjugate points. Hypothesis 4 is commonly called the entropy-gap. Extending the method originally introduced by Gelfert and Ruggiero \cite{gelf19}, we show the following result.
\begin{teo}\label{uniq}
    Let $M$ be a compact $C^{\infty}$ $n$-manifold without conjugate points, $\Tilde{M}$ be its universal covering and $\phi_t$ be its geodesic flow. If $\Tilde{M}$ is a visibility manifold and 
    \[   \sup \{ h_{\mu}(\phi_1): \mu \text{ is supported on } T_1M\setminus \mathcal{R}_0 \}<h(\phi_1), \]
    then the geodesic flow $\phi_t$ has a unique measure of maximal entropy.
\end{teo}
With respect to Theorem \ref{clime}, we retain the entropy-gap, eliminate the third assumption and replace the first two hypothesis by a weaker one: $\tilde{M}$ must be a visibility manifold. Indeed, the existence of a metric of negative curvature implies that $\tilde{M}$ is a visibility manifold by \cite{eber72}, while hypothesis 2 is a necessary consequence of the visibility condition by \cite{rugg03}.

The proof of Theorem \ref{uniq} follows the same reasoning as in Section 7 of \cite{mam23}, which is why we obtained similar results in our setting. We outline the proof and prove the main differences. By Lemma \ref{factor} the geodesic flow $\phi_t$ is time-preserving semi-conjugate to a continuous flow $\psi_t$ under the quotient map $\chi$. By Theorem \ref{frank}, the factor flow $\psi_t$ has a unique measure of maximal entropy $\nu$. We lift $\nu$ to $T_1M$ in order to get the unique measure of maximal entropy for the geodesic flow. Following the measure construction as in \cite{mam23}, let $\epsilon>0$ be an expansivity constant for $\psi_t$. For every $T>0$ we define
$$Per(T,\epsilon)=\{ \chi^{-1}(\gamma): \gamma \text{ is a periodic orbit of }\psi_t \text{ with period in }[T-\epsilon,T+\epsilon]\}.$$
Pick some probability measure $\mu_{\chi^{-1}(\gamma)}$ supported on $\chi^{-1}(\gamma)$ and invariant by the geodesic flow $\phi_t$. We next take the average of the measures
$$\mu_T=\displaystyle\frac{\displaystyle\sum_{\chi^{-1}(\gamma)\in Per(T,\epsilon)}\mu_{\chi^{-1}(\gamma)}}{\# Per(T,\epsilon)}.$$
So, $\mu_T$ is a probability measure on $T_1M$ invariant by the geodesic flow. We consider the measures that are  accumulation points in the weak$^*$ topology of the set $\{ \mu_T:T\to \infty\}$. Following verbatim arguments to \cite{mam23} we get the following property.
\begin{lem}\label{acu}
Let $\mu$ be an accumulation point in the weak$^*$ topology of $\{ \mu_T:T\to \infty\}$ and $\nu$ be the measure of maximal entropy for $\psi_t$. Then, $\mu$ is a measure of maximal entropy for the geodesic flow and $\chi_*\mu=\nu$.
\end{lem}
To show the uniqueness of the measure of maximal entropy we use the following abstract criterion.
\begin{teo}[\cite{buzzi12}]\label{samba}
Let $\phi_t:Y\to Y$ and $\psi_t:X\to X$ be two continuous flows on compact metric spaces, $\chi:Y\to X$ be a time-preserving semi-conjugacy and $\nu$ be the measure of maximal entropy of $\psi_t$. Assume that $\psi_t$ is expansive, has the specification property and  
\begin{enumerate}
\item $h(\phi_1,\chi^{-1}(x))=0$ for every $x\in X$.
\item $\nu \bigg(\{ \chi(y): \chi^{-1}\circ\chi(y)=\{y\} \}\bigg)=1$.
\end{enumerate}
Then $\phi_t$ has a unique measure of maximal entropy.
\end{teo}
To apply this theorem to our case, we only need to show hypothesis 1 and 2. Recall that for every $\eta\in T_1M$, $\chi^{-1}\circ \chi(\eta)=\mathcal{I}(\eta)$. We can thus express hypothesis 1 and 2 in terms of intersections and the expansive set $\mathcal{R}_0$:
$$\sup_{\eta\in T_1M}h(\phi_1,\mathcal{I}(\eta))=0 \quad \text{ and }\quad \nu \bigg(\{ \chi(\eta): \mathcal{I}(\eta)=\{\eta\} \}\bigg)= \nu \bigg( \chi(\mathcal{R}_0)\bigg)=1.$$
The following is an extension of Lemma 6.1 of \cite{gelf20} to our higher dimensional setting and proves hypothesis 1.
\begin{lem}\label{faja}
For every $\eta\in T_1M$, it holds that $h(\phi_1,\mathcal{I}(\eta))=0$.
\end{lem}
\begin{proof}
Let $\tilde{\eta}\in T_1\tilde{M}$ be any lift of $\eta$ hence $\mathcal{\Tilde{I}}(\Tilde{\eta})\subset T_1\Tilde{M}$ is a lift of $\mathcal{I}(\eta)$. We will show that $h(\Tilde{\phi}_1,\mathcal{\Tilde{I}}(\Tilde{\eta}))=0$. This immediately implies the result since $\phi_t$ is a factor of $\Tilde{\phi}_t$. Let $\epsilon>0$, $m\geq 1$ and $E\subset \mathcal{\Tilde{I}}(\Tilde{\eta})$ be a $(m,\epsilon)$-separated set. For each $k=0,\ldots,m-1$ we define the subset $E_k\subset E$ such that if $\eta_1,\eta_2\in E_k$ then $d_s(\Tilde{\phi}_1^k(\eta_1),\Tilde{\phi}_1^k(\eta_2))\geq\epsilon$. Though the sets $E_k$ may have nonempty intersection, we see that $E=\bigcup_{k=0}^{m-1}E_k$. We now estimate the cardinality of $E_k$ for each $k=0,\ldots,m-1$. By Corollary \ref{morse} there exists $\Tilde{Q}>0$ such that $Diam(\Tilde{\phi}_1^k(\mathcal{\Tilde{I}}(\Tilde{\eta})))=Diam(\mathcal{\Tilde{I}}(\Tilde{\phi}_1^k(\Tilde{\eta})))\leq \Tilde{Q}$ hence $\Tilde{\phi}_1^k(\mathcal{\Tilde{I}}(\Tilde{\eta}))$ is contained in a closed ball $\Bar{B}(\xi_k,\Tilde{Q}/2)$ of radius $\Tilde{Q}/2$ centered at some $\xi_k\in T_1\Tilde{M}$. It is known that there exists $A>0$ such that for every $k=0,\ldots,m-1$, $A\epsilon^{-1}\Tilde{Q}^{2n-1}$ is an upper bound for the number of points in $\Bar{B}(\xi_k,\Tilde{Q}/2)$ separated from each other by a distance $\epsilon>0$. So, for every $k=0,\ldots,m-1$, $Card(E_k)\leq A\epsilon^{-1}\Tilde{Q}^{2n-1}$ hence $Card(E)\leq \sum_{k=0}^{m-1}Card(E_k)\leq \sum_{k=0}^{m-1} A\epsilon^{-1}\Tilde{Q}^{2n-1}=mA\epsilon^{-1}\Tilde{Q}^{2n-1}$. Using the entropy definition we get
$$h(\tilde{\phi}_1,\mathcal{\tilde{I}}(\tilde{\eta}))=\lim_{\epsilon\to 0}\limsup_{m\to \infty}\frac{1}{m}\log Card(E)\leq \lim_{\epsilon\to 0}\limsup_{m\to \infty}\frac{1}{m}\log \frac{mA\Tilde{Q}^{2n-1}}{\epsilon}=0.$$
\end{proof}
For hypothesis 2, let us first observe that the geodesic flow $\phi_t$ has positive topological entropy because $\phi_t$ has an expansive factor flow $\psi_t$. Next, we replace Lemma 7.1 of \cite{mam23} by the entropy-gap hypothesis which says that every invariant probability measure supported in $T_1M\setminus \mathcal{R}_0$ has metric entropy $h_{\mu}(\phi_1)$ smaller than topological entropy $h(\phi_1)$. In particular, $T_1M\setminus \mathcal{R}_0$ cannot support a measure of maximal entropy. So, Lemma \ref{acu} implies that any accumulation point $\mu$ of $\{ \mu_T:T\to \infty\}$ is supported in $\mathcal{R}_0$ and hence $\nu=\chi_*\mu$ is supported in $\chi(\mathcal{R}_0)$, which proves hypothesis 2. Thus, applying Theorem \ref{samba} we complete the proof of Theorem \ref{uniq}.    

\section{Visibility manifolds with continuous Green bundles}
This section and the forthcoming ones form the second part of the article where we consider compact manifolds without conjugate points with visibility universal covering, satisfying two additional hypothesis: the continuity of Green bundles and the existence of a hyperbolic periodic point. We will show how these additional hypothesis improve the results of the first six sections. Starting from the abundance of the set where all Lyapunov exponents are non-zero and the hyperbolic periodic points, the regularity of the horospherical foliations and the quotient space, and the uniqueness of the measure of maximal entropy.

In 1977, Pesin \cite{pesin77} established a theory that deals with nonzero Lyapunov exponents. He defined an important set in his theory, which has been called Pesin set ever since. We study a set that is closely related to Pesin set: the set of points where all Lyapunov exponents are nonzero in some subspace transverse to the geodesic flow. We show that this set agrees almost everywhere with an open dense set with respect to Liouville measure. In 1985, Ballmann-Brin-Eberlein \cite{ball85} proved this property for compact rank-1 manifolds of non-positive curvature. In 2003, Ruggiero and Rosas \cite{rosas03} generalize this conclusion to compact manifolds without conjugate points, with bounded asymptote and expansive geodesic flow. Moreover, they showed that periodic hyperbolic points are dense on the unit tangent bundle. We extend these properties to our more general context.

In what follows hyperbolicity and periodicity of points are always referred to the geodesic flow. First, we will prove some basic properties of periodic hyperbolic points in our setting. Recall that for every $\theta\in T_1M$, $G^s(\theta)$ and $G^u(\theta)$ are the stable and unstable Green bundles at $\theta$. Furthermore, $\mathcal{R}_1$ is the set of points where Green bundles are transverse.
\begin{lem}\label{cbas}
Let $M$ be a compact manifold without conjugate points, with visibility universal covering and with continuous Green bundles. If $\eta\in T_1M$ is a periodic hyperbolic point then
\begin{enumerate}
    \item $\eta\in \mathcal{R}_0\cap\mathcal{R}_1$.
    \item $\mathcal{R}_1$ is a non-empty, open and dense set.
    \item If $\Tilde{\eta}\in T_1\Tilde{M}$ is some lift of $\eta$, then for every $\Tilde{\xi}\in \mathcal{\tilde{F}}^s(\Tilde{\eta})\setminus \Tilde{\eta}$ it holds that $d_s(\Tilde{\phi}_t(\Tilde{\eta}),\Tilde{\phi}_t(\Tilde{\xi}))\to \infty$ as $t\to-\infty$.
\end{enumerate}
\end{lem}
\begin{proof}
For item 1, by the hyperbolic structure at $\eta$, $\mathcal{F}^s(\eta)$ and $\mathcal{F}^u(\eta)$ agree locally with the stable and unstable invariant manifolds at $\eta$. Moreover, $G^s(\eta)$ and $G^u(\eta)$ agree with the stable and unstable invariant subspaces at $\eta$, hence they are tangent to $\mathcal{F}^s(\eta)$ and $\mathcal{F}^u(\eta)$ at $\eta$ respectively. Since the invariant subspaces are transverse at $\eta$, so are $G^s(\eta)$ and $G^u(\eta)$, which yields $\eta \in \mathcal{R}_1$. In addition, the above tangency implies that $\mathcal{F}^{s}(\eta)$ and $\mathcal{F}^{u}(\eta)$ are also transverse at $\eta$, which means that $\eta\in \mathcal{R}_0$.

For item 2, item 1 shows that $\mathcal{R}_1$ is non-empty. In addition, the continuity of Green bundles implies that $\mathcal{R}_1$ is an open set. Now, let $U$ be an open set of $T_1M$. By Theorem \ref{v1}(2), we know that $\phi_t$ is topologically transitive hence there exists an orbit $\beta$ of $\phi_t$ that is dense in $T_1M$. Thus $\beta$ intersects both open sets $U$ and $\mathcal{R}_1$. Furthermore, we have $\beta\subset \mathcal{R}_1$ since $\mathcal{R}_1$ is invariant by $\phi_t$. So, $\mathcal{R}_1$ intersects $U$ and therefore $\mathcal{R}_1$ is dense.

For item 3, we suppose by contradiction that there exist $\Tilde{\xi}\in \mathcal{\Tilde{F}}^s(\Tilde{\eta})\setminus \Tilde{\eta}$ such that $d_s(\Tilde{\phi}_t(\Tilde{\eta}),\Tilde{\phi}_t(\Tilde{\xi}))\leq C$ for every $t\leq 0$ and some $C>0$. This implies that $\Tilde{\xi}\in \mathcal{\tilde{F}}^{cu}(\Tilde{\eta})$ and so $\Tilde{\xi}\in \mathcal{\tilde{F}}^s(\Tilde{\eta})\cap \mathcal{\tilde{F}}^{cu}(\Tilde{\eta})$. Using Corollary \ref{key}, we find that $\Tilde{\xi}\in \mathcal{\Tilde{I}}(\Tilde{\eta})$ and hence by projection on $T_1M$ we conclude that $\xi\in \mathcal{I}(\eta)$. By item 1, we know that $\eta$ is expansive, so $\mathcal{I}(\eta)=\{ \eta \}$ and $\xi=\eta$. Seeing that $\Tilde{\xi}\in \mathcal{\Tilde{I}}(\Tilde{\eta})$ we arrive to $\Tilde{\xi}=\Tilde{\eta}$, a contradiction.
\end{proof}
We now define more formally the set mentioned in the introduction,
$$\Lambda_{\phi}=\{ \eta\in T_1M: \chi(\eta,\xi)\neq 0 \text{ for every }\xi\in S_{\eta}\subset T_{\eta}(T_1M) \},$$
where $\chi(\eta,\xi)$ is the Lyapunov exponent of vector $\xi$ at point $\eta$ and $S_{\eta}$ is some subspace transverse to the geodesic flow at $\eta$. For the following results we rely on a Ruggiero's Theorem that we restate below. Recall that $m$ is the Liouville measure defined on $T_1M$.
\begin{teo}[Theorem 4.1 of \cite{rugg21}]\label{t1}
Let $M$ be a compact manifold without conjugate points and with continuous Green bundles. If $\mathcal{R}_1$ is nonempty then for $m$-almost every $\theta \in \mathcal{R}_1$ and every $\xi^s\in G^s(\theta),\xi^u\in G^u(\theta)$, 
$$\lim_{t\to \infty}\frac{1}{t}\log \|d_{\theta}\phi_t(\xi^s)\|<0 \quad \text{ and }\quad  \lim_{t\to \infty}\frac{1}{t}\log \|d_{\theta}\phi_t(\xi^u)\|>0.$$
\end{teo}
\begin{cor}\label{lab}
Let $M$ be a compact manifold without conjugate points, with visibility universal covering and with continuous Green bundles. Suppose the geodesic flow has a periodic hyperbolic point then $\Lambda_{\phi}$ agrees $m$-almost everywhere with the open dense set $\mathcal{R}_1$.
\end{cor}
\begin{proof}
First, observe that by Lemma \ref{cbas}(2), $\mathcal{R}_1$ is a non-empty, open and dense set. So, applying Theorem \ref{t1} to our case we see that for $m$-almost every $\theta\in \mathcal{R}_1$ and every $\xi\in G^s(\theta)\oplus G^u(\theta)$, $\chi(\theta,\xi)\neq 0$. The result follows noting that at any $\theta\in \mathcal{R}_1$, $G^s(\theta)\oplus G^u(\theta)$ span a subspace transverse to the geodesic flow.  
\end{proof}
For another consequence we will need a particular version of a classical Katok's result \cite{katok80}.
\begin{teo}\label{t2}
Let $M$ be a compact manifold, $f:M\to M$ be a $C^2$ diffeomorphism and $\mu$ be a $f$-invariant ergodic measure on $M$. If $\mu$ is not concentrated on a single periodic orbit and $\mu$ has nonzero Lyapunov exponents then for every $x\in Supp(\mu)$ it holds that every open neighborhood of $x$ always contains a periodic hyperbolic point.
\end{teo}
Observe that this theorem is intended for discrete systems. However, there is a large consensus among specialists that Katok Theorem extends to flows without singularities via local cross sections. Furthermore, there are several works that use this extension. Among these works, we can cite to Paternain \cite{pater97}, Barbosa-Ruggiero \cite{barbo13}, Ruggiero \cite{rugg21}, Ledrappier-Lima-Sarig \cite{ledra16} and Araujo-Lima-Poletti \cite{arau20}. 
\begin{cor}
Let $M$ be a compact manifold without conjugate points, with visibility universal covering and with continuous Green bundles. Suppose the geodesic flow has a periodic hyperbolic point then periodic hyperbolic points are dense on $T_1M$.
\end{cor}
\begin{proof}
Since $\mathcal{R}_1$ is invariant by the geodesic flow $\phi_t$, we can restrict the Liouville measure $m$ to $\mathcal{R}_1$. Denote by $m'$ the Liouville measure restricted and normalized to $\mathcal{R}_1$. It follows that $m'$ is a Borel probability measure on $\mathcal{R}_1$ which is invariant by $\phi_t$. Corollary \ref{lab} says that except for the direction tangent to the geodesic flow, $m'$ has nonzero Lyapunov exponents at $m'$-almost every $\theta\in\mathcal{R}_1$. Notice that $m'$ is not concentrated on a single orbit since $\mathcal{R}_1$ is open. Recalling that Pesin \cite{pesin77} essentially proved the ergodicity of Liouville measure when restricted to $\mathcal{R}_1$, we conclude that $m'$ is ergodic. So, the result follows applying Katok Theorem \ref{t2} and noting that $Supp(m')=\mathcal{R}_1$ $m'$-almost everywhere.
\end{proof}

\section{Horospherical foliations are dynamically defined and tangent to Green bundles}\label{s32}
The goal of the section is to study the regularity of the horospherical foliations, the unique integrability of the Green bundles and their tangency to the horospherical foliations.

For closed manifolds of negative curvature, Anosov's work \cite{anos67} implies the uniqueness of continuous foliations of $T_1M$ invariant by geodesic flow. In 1993, Paternain \cite{pater93} showed this uniqueness result for compact surfaces with expansive geodesic flow. This conclusion was extended to higher dimensions under the same hypothesis by Ruggiero in 1997 \cite{rugg97}. In 2007, Barbosa and Ruggiero \cite{barbo07} eliminated the expansivity condition but they still assumed a compact higher genus surface without conjugate points. We extend Barbosa-Ruggiero's Theorem to our higher dimensional setting, which generalizes all previous works. 
\begin{teo}\label{fol}
Let $M$ be a compact $C^{\infty}$ $n$-manifold without conjugate points, with visibility universal covering and with continuous Green bundles. Suppose the geodesic flow has a periodic hyperbolic point $\theta\in T_1M$ then the horospherical foliations $\mathcal{F}^s,\mathcal{F}^u$ are the only continuous invariant $(n-1)$-dimensional foliations with $C^1$-leaves, which are transverse to $\mathcal{F}^u$ (or $\mathcal{F}^s$) at $\theta$.
\end{teo}
Let us prove some preliminary lemmas. Let $\mathcal{D}$ be a foliation satisfying the hypothesis of Theorem \ref{fol}. We first verify the local coincidence of special leaves of $\mathcal{D}$ and $\mathcal{F}^s$ (or $\mathcal{F}^u$). For this we establish as a convention that intersection notations like $\mathcal{D}(\theta)\cap U$ always refer to the connected component of $\mathcal{D}(\theta)\cap U$ containing $\theta$.
\begin{lem}\label{fir}
Let $\theta\in T_1M$ be the periodic hyperbolic point of period $P>0$ given by Theorem \ref{fol} and $\mathcal{D}$ be a $\phi_t$-invariant continuous foliation of $T_1M$ with $C^1$-leaves of dimension $n-1$ which is transverse to $\mathcal{F}^u$ (or $\mathcal{F}^s$) at $\theta$. Then, there exists an open set $U\subset T_1M$ containing $\theta$ such that either $\mathcal{D}(\theta)\cap U=\mathcal{F}^s(\theta)\cap U$ or $\mathcal{D}(\theta)\cap U=\mathcal{F}^u(\theta)\cap U$.
\end{lem}
\begin{proof}
Since $\theta\in T_1M$ is periodic hyperbolic, we choose $U$ as the open neighborhood of $\theta$ where the inclination Lemma holds. Now, without loss of generality, we assume that $\mathcal{D}$ is transverse to $\mathcal{F}^u$ at $\theta$. By contradiction, suppose that 
\begin{equation}\label{aux}
\mathcal{D}(\theta)\cap U\neq \mathcal{F}^s(\theta)\cap U.   
\end{equation}
Let us consider the sets $D_n=\phi_{nP}(\mathcal{D}(\theta))\cap U$ for every $n\geq 1$. Since $\mathcal{D}$ is invariant by the geodesic flow $\phi_t$ we get $\phi_{nP}(\mathcal{D}(\theta))=\mathcal{D}(\theta)$ and hence $D_n\subset \mathcal{D}(\theta) \cap U$ for every $n\geq 1$. From this we conclude that
\begin{equation}\label{cris}
d_r(D_n,\mathcal{F}^u(\theta)\cap U)>0 \quad \text{ for every }n\geq 1, 
\end{equation}
where $d_r$ is the $C^r$-distance for every $r\geq 1$. Now, let us observe that $\theta$ is a hyperbolic fixed point of $\phi_P$ and $\mathcal{D}(\theta)\cap U$ is transverse to $\mathcal{F}^u(\theta)\cap U$ at $\theta$. This together with Equation \eqref{aux} allows us to apply the inclination Lemma to get $d_r(D_n,\mathcal{F}^u(\theta)\cap U)\to 0$ as $n\to\infty$. This contradicts Equation \eqref{cris} and shows that $\mathcal{D}(\theta)\cap U= \mathcal{F}^s(\theta)\cap U$. The other case that assumes that $\mathcal{D}$ is transverse to $\mathcal{F}^s$ at $\theta$ is analogous.
\end{proof}
Let $U$ be the open neighborhood of the periodic hyperbolic point $\theta\in T_1M$ given by this lemma. Consider the following sets
\begin{equation}\label{sets}
    V^s(\theta)= \mathcal{F}^s(\theta)\cap U \quad \text{ and }\quad V^u(\theta)= \mathcal{F}^u(\theta)\cap U.
\end{equation}
We want to generate $\mathcal{F}^s(\theta)$ and $\mathcal{F}^u(\theta)$ with discrete $\phi_t$ iterates of $V^s(\theta)$ and $V^u(\theta)$ respectively. To do this, we work in the covering space $T_1\tilde{M}$. Choosing some lift $\tilde{\theta}\in T_1\tilde{M}$ of $\theta$, we see that $\mathcal{\tilde{F}}^s(\tilde{\theta}),\mathcal{\tilde{F}}^u(\tilde{\theta})$ are lifts of $\mathcal{F}^s(\theta),\mathcal{F}^u(\theta)$ and $\tilde{V}^s(\tilde{\theta}),\tilde{V}^u(\tilde{\theta})$ are also lifts of $V^s(\theta),V^u(\theta)$ such that $\tilde{V}^s(\tilde{\theta})\subset \mathcal{\tilde{F}}^s(\tilde{\theta})$ and $\tilde{V}^u(\tilde{\theta})\subset \mathcal{\tilde{F}}^u(\tilde{\theta})$.

Now, let us note that since $\theta$ is periodic of period $P>0$, we see that $\gamma_{\theta}\subset M$ is a closed geodesic with same period. A basic result in Riemannian geometry \cite{doca92} tell us that $\gamma_{\theta}$ has an associated axial isometry $T:\tilde{M}\to \tilde{M}$ with axis $\gamma_{\tilde{\theta}}$ for some lift $\Tilde{\theta}\in T_1\Tilde{M}$ of $\theta$. This means that $T$ is a translation along $\gamma_{\tilde{\theta}}$: $\gamma_{\tilde{\theta}}(t+P)=T \circ \gamma_{\tilde{\theta}}(t)$ for every $t\in \mathbb{R}$. So, the covering isometry $dT:T_1\tilde{M} \to T_1\tilde{M}$ is also a translation along the orbit of $\tilde{\theta}$: for every $t\in \mathbb{R}$,
\begin{equation}\label{tras}
    \phi_{t+P}(\tilde{\theta})= dT\circ \phi_t(\tilde{\theta}).
\end{equation}
\begin{lem}
Let $\theta\in T_1M$ be the periodic hyperbolic point of period $P>0$ given by Theorem \ref{fol}, $\tilde{\theta}\in T_1\Tilde{M}$ be any lift of $\theta$, $T$ be the axial isometry associated to $\gamma_{\theta}$ having axis $\gamma_{\tilde{\theta}}$ and $\tilde{V}^s(\tilde{\theta}),\tilde{V}^u(\tilde{\theta})$ be the lifts defined above. Then, $(dT)^{-n}\circ \phi_{nP}(\tilde{\theta})=\tilde{\theta}$ for every $n\in \mathbb{Z}$ and 
\begin{equation}\label{dentro}
    \mathcal{\tilde{F}}^s(\tilde{\theta})=\bigcup_{n\in \mathbb{N}}(dT)^n\circ \phi_{-nP}(\tilde{V}^s(\tilde{\theta})) \quad \text{ and }\quad \mathcal{\tilde{F}}^u(\tilde{\theta})=\bigcup_{n\in \mathbb{N}}(dT)^{-n}\circ \phi_{nP}(\tilde{V}^u(\tilde{\theta})).
    \end{equation}
\end{lem}
\begin{proof}
For the first assertion, Equation \eqref{tras} shows that for every $t\in \mathbb{R}$ and every $n\in \mathbb{Z}$,
$\phi_{t+nP}(\tilde{\theta})=(dT)^{n}\circ \phi_t(\tilde{\theta})$ and hence $(dT)^{-n}\circ \phi_{t+nP}(\tilde{\theta})=\phi_t(\tilde{\theta})$. The conclusion follows by setting $t=0$. To prove the direct inclusion of the first expression in Equation \eqref{dentro}, note that horospherical foliations are invariant by the geodesic flow $\phi_t$ and by all covering isometries of $T_1\Tilde{M}$. So, for every $n\geq 1$ we have 
$$(dT)^{n}\circ \phi_{-nP}(\mathcal{\tilde{F}}^s(\tilde{\theta}))= \mathcal{\tilde{F}}^s((dT)^{n}\circ \phi_{-nP}(\tilde{\theta}))=\mathcal{\tilde{F}}^s(\tilde{\theta}).$$
Since $\tilde{V}^s(\tilde{\theta})\subset \mathcal{\tilde{F}}^s(\tilde{\theta})$, we see that for every $n\geq 1$, $(dT)^{n}\circ \phi_{-nP}(\tilde{V}^s(\tilde{\theta}))\subset \mathcal{\tilde{F}}^s(\tilde{\theta})$ and hence  $\bigcup_{n\in \mathbb{N}}(dT)^n\circ \phi_{-nP}(\tilde{V}^s(\tilde{\theta}))\subset \mathcal{\tilde{F}}^s(\tilde{\theta})$.
For the reverse inclusion, Lemma \ref{cbas}(3) states that for every $\tilde{\eta}\in \tilde{V}^s(\tilde{\theta})\subset \mathcal{\tilde{F}}^s(\tilde{\theta})$, $d_s(\phi_t(\tilde{\theta}),\phi_t(\tilde{\eta}))\to \infty$ as $t\to -\infty$.
Thus, setting $t=-nP$ for every $n\geq 1$ we get
\begin{align*}
  d_s(\tilde{\theta},(dT)^{n}\circ \phi_{-nP}(\tilde{\eta}))&=d_s((dT)^{n}\circ \phi_{-nP}(\tilde{\theta}),(dT)^{n}\circ \phi_{-nP}(\tilde{\eta}))\\
  &=d_s( \phi_{-nP}(\tilde{\theta}), \phi_{-nP}(\tilde{\eta}))\to \infty \quad\text{ as }\quad n\to \infty.  
\end{align*}
This implies that $Diam((dT)^n\circ \phi_{-nP}(\tilde{V}^s(\tilde{\theta})))\to \infty$ as $n\to \infty$. From this we conclude that $\mathcal{\tilde{F}}^s(\tilde{\theta})$ is covered by $\bigcup_{n\in \mathbb{N}}(dT)^n\circ \phi_{-nP}(\tilde{V}^s(\tilde{\theta}))$ as $n\to \infty$, which proves the inclusion. The proof for the other expression is analogous.
\end{proof}
As an immediate application the entire leaves $\mathcal{D}(\theta)$ and $\mathcal{F}^s(\theta)$ agree for the periodic hyperbolic point $\theta$.
\begin{lem}\label{coin}
Let $\theta\in T_1M$ be the periodic hyperbolic point of period $P>0$ given by Theorem \ref{fol} and $\mathcal{D}$ be a $\phi_t$-invariant continuous foliation of $T_1M$ with $C^1$-leaves of dimension $n-1$ which is transverse to $\mathcal{F}^u$ (or $\mathcal{F}^s$) at $\theta$. Then,

\begin{enumerate}
    \item If $V^s(\theta)$ and $V^u(\theta)$ are the sets defined in Equation \eqref{sets} then 
    $$\mathcal{F}^s(\theta) =\bigcup_{n\in \mathbb{N}}\phi_{-nP}(V^s(\theta)) \quad\text{ and }\quad  \mathcal{F}^u(\theta)=\bigcup_{n\in \mathbb{N}}\phi_{nP}(V^u(\theta)).$$
    \item Either $\mathcal{D}(\theta)=\mathcal{F}^s(\theta)$ or $\mathcal{D}(\theta)=\mathcal{F}^u(\theta)$.
\end{enumerate}
\end{lem}
\begin{proof}
Item 1 follows by applying the covering map $d\pi$ to Equation \eqref{dentro} since
$$d\pi(\bigcup_{n\in \mathbb{N}}(dT)^n\circ \phi_{-nP}(\tilde{V}^s(\tilde{\theta})))=\bigcup_{n\in \mathbb{N}}d\pi \circ (dT)^n\circ \phi_{-nP}(\tilde{V}^s(\tilde{\theta}))=\bigcup_{n\in \mathbb{N}} \phi_{-nP}(\tilde{V}^s(\tilde{\theta})).$$
For item 2, consider any $\eta \in \mathcal{F}^s(\theta)$. By item 1 we know that $\eta\in \phi_{-nP}(V^s(\theta))$ for some $n\geq 1$. Without loss of generality, by Lemma \ref{fir} we suppose that $\mathcal{D}(\theta)=\mathcal{F}^s(\theta)$ on $V^s(\theta)$. Since $\mathcal{D}$ and $\mathcal{F}^s$ are invariant by $\phi_{-nP}$ for every $n\geq 1$, it follows that $\mathcal{D}(\theta)=\mathcal{F}^s(\theta)$ on $\phi_{-nP}(V^s(\theta))$. So, by item 1 $\mathcal{D}(\theta)=\mathcal{F}^s(\theta)$ agree on $\phi_{-NP}(V^s(\theta))$ containing $\eta$ for some $N\geq 1$. Therefore $\mathcal{D}(\theta)=\mathcal{F}^s(\theta)$ since $\eta$ was arbitrary. The other assumption $\mathcal{D}(\theta)=\mathcal{F}^u(\theta)$ on $V^u(\theta)$ is analogous.
\end{proof}
Now, using the density of horospherical leaves (Theorem \ref{v1}(1)) we can extend item 2 to any arbitrary $\eta\in T_1M$.
\begin{proof}[Proof of Theorem \ref{fol}]
Let $\mathcal{D}$ be a $\phi_t$-invariant continuous foliation of $T_1M$ with $C^1$-leaves of dimension $n-1$, $\xi\in T_1M$ and $B$ be a closed ball containing $\xi$. Without loss of generality, we assume that $\mathcal{D}$ is transverse to $\mathcal{F}^u$ at the hyperbolic periodic point $\theta\in T_1M$. Applying Lemma \ref{coin} we obtain
\begin{equation}\label{coin1}
  \mathcal{D}(\theta)=\mathcal{F}^s(\theta).  
\end{equation}
By Theorem \ref{v1}(1) we know that $\mathcal{F}^s(\theta)$ is dense on $T_1M$ hence there exists a sequence $\xi_n\in \mathcal{F}^s(\theta)$ converging to $\xi$. For every $n\geq 1$, let $F_n$ and $D_n$ be the respective connected components of $\mathcal{F}^s(\xi_n)\cap B$ and $\mathcal{D}(\xi_n)\cap B$ containing $\xi_n$. Equation \eqref{coin1} provides that $F_n=D_n$ for every $n\geq 1$. Denote by $F$ and $D$ the respective connected components of $\mathcal{F}^s(\xi)\cap B$ and $\mathcal{D}(\xi)\cap B$ containing $\xi$. Since $\mathcal{F}^s$ and $\mathcal{D}$ are continuous foliations in the compact-open topology, we deduce that $F_n$ converges to $F$ and $F_n=D_n$ converges to $D$ as $n\to \infty$ in the compact-open topology. This implies that $F=D$ and hence $\mathcal{F}^s(\xi)\cap B=\mathcal{D}(\xi)\cap B$ for every $\xi\in T_1M$.
This means that foliations $\mathcal{D}$ and $\mathcal{F}^s$ agree. A similar reasoning shows that $\mathcal{D}$ and $\mathcal{F}^u$ agree if we assume that $\mathcal{D}(\theta)=\mathcal{F}^u(\theta)$. 
\end{proof}
We now turn to the problem of the tangency between horospherical foliations and Green bundles. We mention that in the context of compact manifolds without conjugate points, Knieper \cite{knip86} studied a related problem: the integrability of Green bundles. Assuming furthermore the continuity of Green bundles, Knieper showed that Green bundles integrate to continuous invariant foliations. Moreover, he conjectured the tangency between horospherical foliations and Green bundles under these hypothesis.
\begin{teo}[\cite{knip86}]\label{greni}
Let $M$ be a compact $n$-manifold without conjugate points. If Green bundles are continuous then these bundles integrate to continuous invariant $(n-1)$-dimensional foliations $\mathcal{G}^s$ and $\mathcal{G}^u$ of $T_1M$ with $C^1$-leaves.
\end{teo}
Note that this theorem does not ensure the unique integrability of Green bundles. An application of Theorem \ref{fol} allows us to answer both problems: unique integrability and tangency. 
\begin{cor}\label{tang}
Let $M$ be a compact $C^{\infty}$ $n$-manifold without conjugate points, with visibility universal covering and with continuous Green bundles. Suppose the geodesic flow has a periodic hyperbolic point then Green bundles are uniquely integrable and tangent to the horospherical foliations. In particular, $\mathcal{R}_1\subset \mathcal{R}_0$.
\end{cor}
\begin{proof}
By Theorem \ref{greni}, Green bundles integrate to some $\phi_t$-invariant continuous foliations $\mathcal{G}^s$ and $\mathcal{G}^u$ of $T_1M$ with $C^1$-leaves of dimension $n-1$. By hypothesis there exists a hyperbolic periodic point $\theta\in T_1M$ hence by the hyperbolic structure at $\theta$,
\begin{equation}\label{hyp1}
  \mathcal{G}^s(\theta)\cap U=\mathcal{F}^s(\theta)\cap U\quad \text{ and }\quad  \mathcal{G}^u(\theta)\cap U=\mathcal{F}^u(\theta)\cap U,
\end{equation}
for some open set $U\subset T_1M$ containing $\theta$. From Lemma \ref{cbas}(1) we know that $\mathcal{G}^s$ ($\mathcal{G}^u$) is transverse to $\mathcal{F}^u$ ($\mathcal{F}^s$) at $\theta$. So applying Theorem \ref{fol} to $\mathcal{G}^s$, either $\mathcal{G}^s=\mathcal{F}^s$ or $\mathcal{G}^s=\mathcal{F}^u$. But Equation \eqref{hyp1} ensures that $\mathcal{G}^s=\mathcal{F}^s$ and $\mathcal{G}^u=\mathcal{F}^u$. The conclusion follows since $\mathcal{G}^s$ and $\mathcal{G}^u$ are arbitrary foliations that were integrated from Green bundles.
\end{proof}

\section{A topological manifold structure for the quotient space}
In this section we study the regularity of the quotient space defined in Section \ref{c26} under the assumption of the continuity of Green bundles and the existence of a hyperbolic periodic orbit.

Recall that $\psi_t$ is a continuous flow acting on the quotient space $X$ and $\chi:T_1M \to X$ is the corresponding quotient map. By Lemma \ref{factor}, $X$ is in general a compact metric space of topological dimension at least $n-1$. The topological manifold structure of $X$ was first proved by Gelfert and Ruggiero for compact higher genus surfaces without focal points and for compact higher genus surfaces without conjugate points and with continuous Green bundles \cite{gelf19,gelf20}. We extend this result to our higher dimensional context.
\begin{teo}\label{top}
Let $M$ be a compact $C^{\infty}$ $n$-manifold without conjugate points, with visibility universal covering and with continuous Green bundles. Suppose the geodesic flow has a periodic hyperbolic point then quotient space $X$ is a $(2n-1)$-dimensional compact topological manifold.    
\end{teo}
We first show the existence of special stable and unstable horospherical leaves.
\begin{lem}\label{inc}
Let $M$ be a compact $C^{\infty}$ $n$-manifold without conjugate points, with visibility universal covering and with continuous Green bundles. If $\theta\in T_1M$ is a periodic hyperbolic point of period $P>0$ then horospherical leaves $\mathcal{F}^s(\theta)$ and $\mathcal{F}^u(\theta)$ are included in the expansive set $\mathcal{R}_0$.
\end{lem}
\begin{proof}
Lemma \ref{cbas} says that $\theta\in \mathcal{R}_1$ and moreover guarantees the existence of an open neighborhood $U$ of $\theta$ totally contained in $\mathcal{R}_1$. Repeating the same argument of Lemma \ref{coin}(1) with $U$ instead of $V(\theta)$, we have $\mathcal{F}^s(\theta) =\bigcup_{n\in \mathbb{N}}\phi_{-nP}(U)$. Since $\mathcal{R}_1$ is invariant by the geodesic flow, we obtain $\mathcal{F}^s(\theta)\subset \mathcal{R}_1$. Finally, Corollary \ref{tang} provides $\mathcal{F}^s(\theta)\subset\mathcal{R}_0$. The unstable case is analogous.
\end{proof}

\begin{figure}[ht]
    \centering
    \includegraphics[width=0.6\textwidth]{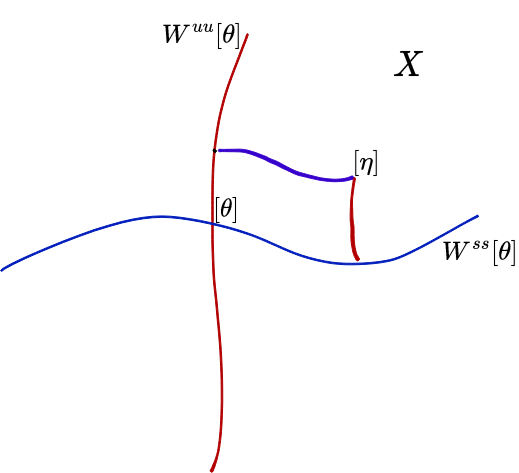}
    \caption{Parametrization by the local product neighborhood of $[\theta]$ in the quotient space $X$.}
    \label{prod}
\end{figure}

Recall that for every $[\eta]\in X$, the strong stable, the strong unstable, the center stable and center unstable sets of $[\eta]$ are denoted by $W^{ss}[\eta]$, $W^{uu}[\eta]$, $W^{cs}[\eta]$ and $W^{cs}[\eta]$ respectively. From Proposition \ref{loc}(2) we know that $W^{ss}[\eta]$ and $W^{uu}[\eta]$ agree with $\chi(\mathcal{F}^s(\eta))$ and $\chi(\mathcal{F}^u(\eta))$ locally. We will use the local product structure to construct a 'local product neighborhood' and hence a local Euclidean neighborhood of any point (See Figure \ref{prod}). 

\begin{proof}[Proof of Theorem \ref{top}]
Since $X$ is Hausdorff and second-countable it only remains to show that $X$ is locally Euclidean. Let $[\eta]\in X$ and $B([\eta],\epsilon')$ be an open ball. By Theorem \ref{v1}(1), $\mathcal{F}^s(\xi)$ and $\mathcal{F}^u(\xi)$ intersect $\chi^{-1}(B([\eta],\epsilon'))$ hence $W^{ss}[\xi]$ and $W^{uu}[\xi]$ intersect $B([\eta],\epsilon')$ for $\epsilon'$ small enough. Let $[\xi_s]\in W^{ss}[\xi]\cap B([\eta],\epsilon')$ and $[\xi_u]\in W^{uu}[\xi]\cap B([\eta],\epsilon')$. We denote by $W^{s}[\xi_s]$ and $W^{u}[\xi_u]$ the connected components of $W^{ss}[\xi]\cap B([\eta],\epsilon')$ and $W^{uu}[\xi]\cap B([\eta],\epsilon')$ containing $[\xi_s]$ and $[\xi_u]$ respectively. Let $\Tilde{\xi}_s,\Tilde{\xi}_u\in T_1\Tilde{M}$ be some lifts of $\xi_s,\xi_u$. Corollary \ref{par} says that both $\mathcal{\Tilde{F}}^s(\Tilde{\xi}_s)$ and $\mathcal{\Tilde{F}}^u(\Tilde{\xi}_u)$ are locally parametrized by $(n-1)$ polar coordinates each. This means that there exist homeomorphisms $\Tilde{f}_s:B_s(0,\epsilon)\subset \mathbb{R}^{n-1}\to\mathcal{\Tilde{F}}^s(\Tilde{\xi}_s)$ and $\Tilde{f}_u:B_u(0,\epsilon)\subset \mathbb{R}^{n-1}\to\mathcal{\Tilde{F}}^u(\Tilde{\xi}_u)$ onto their images where $B_s(0,\epsilon)$ and $B_u(0,\epsilon)$ are open balls. Choosing $\epsilon>0$ small enough, the maps $d\pi\circ \Tilde{f}_s:B_s(0,\epsilon) \to \mathcal{F}^s(\xi_s)$ and $d\pi\circ \Tilde{f}_u:B_u(0,\epsilon)\to \mathcal{F}^u(\xi_u)$ are homeomorphisms onto their images because $d\pi$ is a local homeomorphism. Since $\mathcal{F}^s(\xi_s),\mathcal{F}^u(\xi_u)\subset \mathcal{R}_0$ by Lemma \ref{inc}, we see that quotient map $\chi$ is bijective hence a homeomorphism when restricted to $\mathcal{F}^s(\xi_s)$ and $\mathcal{F}^u(\xi_u)$. So, the maps $f_s=\chi\circ d\pi\circ \Tilde{f}_s:B_s(0,\epsilon)\subset \mathbb{R}^{n-1}\to U_s\subset W^s[\xi_s]$ and $f_u=\chi\circ d\pi\circ \Tilde{f}_u:B_u(0,\epsilon)\subset \mathbb{R}^{n-1}\to U_u\subset W^u[\xi_u]$    
are homeomorphisms onto their images. For $\tau>0$ small enough, we can define 
\[ f:B_s(0,\epsilon)\times B_u(0,\epsilon)\times (-\tau,\tau) \to U\subset X, \; (x,y,t)\mapsto \psi_t(W^s(f_s(x))\cap W^{cu}(f_u(x))).  \] 
Since $f_s,f_u$ and the local product structure of $\psi_t$ are continuous, we deduce that $f$ is continuous as well. Thus, $f$ is a continuous parametrization of $U\subset X$ by an open set of $\mathbb{R}^{2n-1}$. By the invariant domain Theorem we conclude that $U$ is an open set of $X$. We can also recover the 'coordinates' of every $[\theta]\in U$ with the following formulas:
\[ [\theta_s]=W^s[\xi_s]\cap W^{cu}[\theta]\in W^s[\xi_s], \qquad [\theta_u]=W^{cs}[\theta]\cap W^{u}[\xi_u]\in W^{u}[\xi_u] \] 
and $t\in (-\tau,\tau)$ such that $[\theta_u]=W^{s}(\psi_t[\theta])\cap W^{u}[\xi_u]$. Therefore, the map
\[ g:U\subset X\to B_s(0,\epsilon)\times B_u(0,\epsilon)\times (-\tau,\tau)\subset \mathbb{R}^{2n-1}, \quad [\theta]\mapsto (f_s^{-1}[\theta_s],f_u^{-1}[\theta_u],t)\]
is an inverse for $f$. By the continuity of the local product structure we similarly find that $g$ is continuous and so $f$ is a homeomorphism. We note that choosing $\epsilon'>0$ small enough we can pick $W^s[\xi_s]$ and $W^u[\xi_u]$ such that $[\eta_s]\in U_s$ and $[\eta_u]\in U_u$ hence $[\eta]\in U$. Thus, $U\subset X$ is an open neighborhood of $[\eta]$ that is homeomorphic to an open set of $\mathbb{R}^{2n-1}$ hence $X$ is locally Euclidean.
\end{proof}
We remark that continuity of Green bundles allowed us to construct the homeomorphisms from open sets of $\mathbb{R}^{n-1}$ into some strong stable and strong unstable sets of the quotient flow: 'the coordinates axes'. 

\section{Uniqueness of measure of maximal entropy for some manifolds without conjugate points}
This section is devoted to prove that certain family of manifolds satisfy the entropy-gap hypothesis and hence they have the following property.
\begin{teo}\label{eje}
Let $M$ be a compact $C^{\infty}$ $n$-manifold without conjugate points, with visibility universal covering and with continuous Green bundles. Suppose the geodesic flow $\phi_t$ has a periodic hyperbolic point then $\phi_t$ has a unique measure of maximal entropy which has full support.
\end{teo}
This theorem gives a family of manifolds with unique measures of maximal entropy which includes compact manifolds without focal points and compact manifolds of bounded asymptote, both admitting a hyperbolic periodic orbit for their geodesic flows.

Recall that under our hypothesis, we have a quotient flow $\psi_t$ acting on a compact topological $(2n-1)$-manifold $X$ and a time-preserving semi-conjugacy $\chi:T_1M\to X$. Moreover, by Theorem \ref{frank} $\psi_t$ has a unique measure of maximal entropy $\nu$. The following restatement of Proposition 7.3.15 of \cite{fish19} in our context will imply that $\nu$ has full support on $X$.
\begin{pro}\label{gibbs}
Let $X$ be a compact metric space, $\psi_t:X\to X$ be a continuous expansive flow with the specification property and $\nu$ be its unique measure of maximal entropy. Then, for every $\epsilon>0$ there exist $A_{\epsilon},B_{\epsilon}>0$ such that for every $x\in X$ and every $T\geq 0$, it holds $A_{\epsilon}\leq e^{Th(\psi_1)}\nu(B(x,\epsilon,T))\leq B_{\epsilon}$ where $B(x,\epsilon,T)$ is a $(T,\epsilon)$-dynamical ball. 
\end{pro}
\begin{cor}\label{supp}
Let $M$ be a compact $C^{\infty}$ $n$-manifold without conjugate points which has a visibility universal covering and $\psi_t:X\to X$ be quotient flow. Then, the unique measure of maximal entropy $\nu$ for $\psi_t$ has full support.
\end{cor}
\begin{proof}
For $T=0$ and every $\epsilon>0$, by Theorem \ref{prop} we can apply Proposition \ref{gibbs} to $\psi_t$ and $\nu$. We thus get $0<A_{\epsilon}\leq \nu(B(x,\epsilon,0))\leq B_{\epsilon}$ for every $\epsilon>0$ and every $x\in X$. Hence $\nu$ has full support on $X$ since $B(x,\epsilon,0)$ is just an open ball of radius $\epsilon$ centered at $x$.    
\end{proof}
We now show that push-forward map $\chi_*$ carries measures of maximal entropy into measures of maximal entropy.
\begin{lem}\label{mea}
Let $M$ be a compact $C^{\infty}$ $n$-manifold without conjugate points which has a visibility universal covering, $\nu$ be the unique measure of maximal entropy for $\psi_t$ and $\chi$ be the quotient map. If $\mu$ is a measure of maximal entropy for the geodesic flow then so is $\chi_*\mu$ for $\psi_t$. In particular, $\chi_*\mu=\nu$ for every $\mu$. 
\end{lem}
\begin{proof}
Clearly $\chi_*\mu$ is a probability Borel measure on $X$ invariant by $\psi_t$. Since $(T_1M,\phi_t)$ is an extension of $(X,\psi_t)$, we have $h(\psi_1)\leq h(\phi_1)$. Similarly, since $(T_1M,\phi_t,\mu)$ is an extension of $(X,\psi_t,\chi_*\mu)$, $h_{\chi*\mu}(\psi_1)\leq h_{\mu}(\phi_1)$. As all intersections $\mathcal{I}(\eta)$ carry no topological entropy by Lemma \ref{faja}, Bowen's formula \cite{bowen71endo} and Ledrappier-Walter's formula \cite{ledra77} imply that
\[ h(\phi_1)\leq h(\psi_1)+\sup_{[\eta]\in X} h(\phi_1,\chi^{-1}[\eta])= h(\psi_1)+\sup_{[\eta]\in X} h(\phi_1,\mathcal{I}(\eta))= h(\psi_1).\]
\[ h_{\mu}(\phi_1)\leq h_{\chi*\mu}(\psi_1)+ \int_X h(\phi_1,\chi^{-1}[\eta])d\chi_*\mu = h_{\chi*\mu}(\psi_1).\]
We thus get $h(\phi_1)=h(\psi_1)$ and $h_{\mu}(\phi_1)=h_{\chi_*\mu}(\psi_1)$. Because $\mu$ is a measure of maximal entropy for $\phi_t$, the calculation $h_{\chi_*\mu}(\psi_1)=h_{\mu}(\phi_1)=h(\phi_1)=h(\psi_1)$ shows that $\chi_*\mu$ is a measure of maximal entropy for $\psi_t$.
\end{proof}
Assuming more hypothesis on the manifold $M$ we can show the following feature of measures of maximal entropy.
\begin{lem}\label{supo}
Under the hypothesis of Theorem \ref{eje}, every measure $\mu$ of maximal entropy for the geodesic flow has full support.
\end{lem}
\begin{proof}
Let $U$ be any open set of $T_1M$. By Lemma \ref{cbas} and Corollary \ref{tang} we know that $\mathcal{R}_1$ is a non-empty dense set included in $\mathcal{R}_0$. So, there is an expansive point $\xi\in U$ and hence there exists an open saturated set $U'\subset U$ by Lemma \ref{sat}(2). Let $\nu$ be the unique measure of maximal entropy for $\psi_t$ and $\chi$ be the quotient map. Since $\chi(U')$ is an open set of $X$, Lemma \ref{mea} and Corollary \ref{supp} complete the proof: $\mu(U)\geq \mu(U')=\mu(\chi^{-1}\circ\chi(U'))=\chi_*\mu(\chi(U'))=\nu(\chi(U'))>0$.
\end{proof}
Recall that $\mathcal{R}_1$ is the set of points where Green bundles are transverse.
\begin{lem}\label{canno}
Under the hypothesis of Theorem \ref{eje}, $T_1M\setminus \mathcal{R}_0$ and $T_1M\setminus \mathcal{R}_1$ cannot support a measure of maximal entropy for the geodesic flow.
\end{lem}
\begin{proof}
Suppose by contradiction that $T_1M\setminus \mathcal{R}_1$ supports a measure of maximal entropy $\mu$. From Lemma \ref{supo} we know that $\mu$ has full support. Since $\mathcal{R}_1$ is a non-empty open set by Lemma \ref{cbas}(2), we have $\mu(\mathcal{R}_1)>0$, a contradiction to the support of $\mu$. The conclusion for $T_1M\setminus \mathcal{R}_0$ follows from $\mathcal{R}_1\subset \mathcal{R}_0$ given by Corollary \ref{tang}.
\end{proof}
We highlight that in our context, the compactness of $T_1M\setminus \mathcal{R}_1$ provides the following variational principle for $T_1M\setminus \mathcal{R}_1$: $h(\phi_1,T_1M\setminus \mathcal{R}_1)=\sup_{\mu}h_{\mu}(\phi_1,T_1M\setminus \mathcal{R}_1)$, where $\mu$ varies over all probability invariant measures supported on $T_1M\setminus \mathcal{R}_1$.
\begin{lem}\label{hipo}
Under the hypothesis of Theorem \ref{eje}, $h(\phi_1,T_1M\setminus\mathcal{R}_1)<h(\phi_1)$. In particular, the entropy-gap assumption of Theorem \ref{uniq} is satisfied.  
\end{lem}
\begin{proof}
Suppose by contradiction that $h(\phi_1,T_1M\setminus \mathcal{R}_1)=h(\phi_1)$. Using the variational principle, we can construct a sequence of invariant measures $\mu_n$ supported on $T_1M\setminus \mathcal{R}_1$ with $h_{\mu_n}(\phi_1,T_1M\setminus \mathcal{R}_1)\to h(\phi_1)$. Since $T_1M\setminus\mathcal{R}_1$ is compact we conclude the set $\mathcal{M}$ of probability invariant measures supported on $T_1M\setminus\mathcal{R}_1$ is weak$^*$ compact. So, there exists a subsequence $\mu_{n_k}$ converging weakly$^*$ to some $\mu\in \mathcal{M}$. Since manifold $M$ is smooth, Newhouse's work \cite{new89} provides that entropy map restricted to $\mathcal{M}$ is upper semi-continuous and hence $h_{\mu}(\phi_1)\geq \lim_k h_{\mu_{n_k}}(\phi_1,T_1M\setminus\mathcal{R}_1)=h(\phi_1)$. Therefore $\mu$ is a measure of maximal entropy supported on $T_1M\setminus \mathcal{R}_1$ which contradicts Lemma \ref{canno}. 
\end{proof}
\begin{proof}[Proof of Theorem \ref{eje}]
Notice that together with Lemma \ref{hipo} all the hypothesis of Theorem \ref{uniq} are satisfied hence the geodesic flow has a unique measure of maximal entropy which has full support by Lemma \ref{supo}.
\end{proof}

\section{Acknowledgments}
The first author appreciates the financial support of CNPQ and FAPEMIG funding agencies during the work. The second author was financed by CNPQ, Edital Universal, Faperj and PRONEX-Geometria. This article was supported in part by INCTMat under the project INCTMat-Faperj (E26/200.866/2018) and by Sorbonne Université for the project Emergence Dyflo 2021-2022.

\bibliographystyle{plain}
\bibliography{ref}

\end{document}